\documentclass{article}

\usepackage{mathpazo}
\usepackage{graphicx}
\usepackage{amsmath,amssymb,amsfonts}
\usepackage{amsthm}
\usepackage{mathrsfs}
\usepackage[title]{appendix}
\usepackage{xcolor}
\usepackage{textcomp}
\usepackage{manyfoot}
\usepackage{booktabs}
\usepackage{algorithm}
\usepackage{algorithmicx}
\usepackage{algpseudocode}
\usepackage[colorlinks]{hyperref}

\usepackage{tikz-cd}
\usepackage{tikz}
\usetikzlibrary{tikzmark}
\usetikzlibrary{arrows}
\usetikzlibrary{decorations.pathreplacing}
\usetikzlibrary{decorations.markings}

\input xy    
\xyoption{curve}  
\xyoption{arc} 
\xyoption{frame} 
\xyoption{tips} 
\xyoption{arrow} 
\xyoption{color} 
\xyoption{poly}

\newtheorem{thm}{Theorem}[section]
\newtheorem{prop}[thm]{Proposition} 
\newtheorem{cor}[thm]{Corollary}
\newtheorem{lemma}[thm]{Lemma}

\newtheorem{ex}{Example}
\newtheorem{remark}{Remark}

\newtheorem*{theorem*}{Theorem}
\newtheorem{theorem}{Theorem}

\newcommand{\defeq}{\mathrel{\mathpalette{\vcenter{\hbox{$:$}}}=}}
\newcommand{\A}{\mathcal{A}}
\newcommand{\G}{\mathcal{G}}
\newcommand{\fix}{\mathrm{Fix}}
\newcommand{\SO}{\mathrm{SO}}

\newcommand{\bb}[1]{\mathbb{#1}}

\newcommand{\sumdd}[2]{\displaystyle \sum_{#1}^{#2}}
\newcommand{\homeo}{\text{Homeo}}
\newcommand{\R}{\mathcal R}
\newcommand{\X}{\mathfrak X}

\tikzset{pics/.cd,
handle/.style={code={
\draw  (-2,0) coordinate (-lpoint) 
to [out=100, in=300] (-2.8,2) 
to [out=120, in=250] (-3,4) 
to [out=70,in=180] (0,6) 
to [out=0,in=110] (3,4) 
to [out=290,in=60] (2.8,2) 
to [out=240,in=80] (2,0)  coordinate (-rpoint);
\draw (0,.75) .. controls (1,2) and (1,3.25) .. (0,4.25);
\draw (0.2,1) .. controls (-.75,2.25) and (-.75,3.25) .. (0.25,4);
\draw[blue,->, line width=3pt] (0,-pi) coordinate  
to [out=100, in=270] (-1.5,3.25) 
to [out=90, in=180] (0,5) ;
\draw[blue, line width=3pt]  (0,5) coordinate  
to [out=0,in=90] (1.5,3.25) 
to [out=270,in=80] (0,-pi);
\draw[red,->, line width=3pt] (0,-pi) coordinate  
to [out=60, in=230] (1.5,-1) ;
\draw[red, line width=3pt] (1.5,-1) coordinate  
to [out=50, in=250] (2,0) ;
\draw[red, loosely dotted,line width=3pt] (2,0) coordinate  
to [out=100, in=80] (.25,1) ;
\draw[red, line width=3pt] (.25,1) coordinate  
to [out=260, in=90] (0,-pi) ;
\filldraw[black] (0,-pi) circle (6pt);
\draw[-to,line width=2pt] (7.5, -2.5) arc
    [
        start angle=-7,
        end angle=75,
        x radius=5cm,
        y radius =5cm
    ] ;
}}}

\title{Semifree actions and representation varieties}

\author{Patrick Dendis, David L. Duncan, Timothy Tarter, \\ Lior Weber, and Colin Willette}

\date{}

\begin{document}

\maketitle

\begin{abstract}
A group action on a topological space induces an action on the space's representation variety. We determine the fixed point sets of this latter action for semifree group actions on surfaces and circle bundles over surfaces. We also discuss implications for mapping tori of periodic maps. 
\end{abstract}

\tableofcontents

\section{Introduction}

The interaction between finite group actions and low-dimensional topology has a long and storied history that includes the resolution of numerous famous problems such as the Smith Conjecture \cite{Bass-Morgan} and spherical space forms problem (a consequence of the {geometrization conjecture} \cite{P1,P2,P3}), as well as many as-yet unresolved problems, such as Problems 3.75, 4.84 and 4.86 in the recent list \cite{K3}. Likewise, representation varieties have been used by low-dimensional topologists for decades. In particular, the $S^3$-representation variety $\R(X) \defeq \hom(\pi_1(X), S^3)$ and character variety $\X(X) \defeq \R(X) / S^3$ of a manifold $X$ have repeatedly proven to be powerful weapons in the low-dimensional topologist's arsenal \cite{AM}, \cite{Floer}, \cite{SZ}, \cite{LP-CZ}. 

Throughout, $\Sigma$ will denote a surface that is closed, connected and oriented, and we assume this is equipped with an action of a finite group $\pi$ that is effective and semifree, and acts through orientation-preserving homeomorphisms. The general aim of this article is to investigate the induced $\pi$-actions on $\R(\Sigma)$ and $\X(\Sigma)$; we primarily focus on determining the topological type of the fixed point sets. Section \ref{sec:InGeneral} explores this problem in general, while the earlier Section \ref{sec:SpecialCases} gives a detailed analysis of various interesting special cases. For each of these special cases, we also describe the $\pi$-fixed point sets in $\R(N)$ and $\X(N)$ for any $\pi$-equivariant oriented circle bundle $N$ over $\Sigma$. The fixed point sets on $\X(\Sigma)$ are closely related to the representation and character varieties of the mapping torus $\Sigma_\phi$ for any generator $\phi \in \pi$, and we describe these as well. As we will see, the components of the various fixed point sets tend to be familiar spaces such as points, spheres, or intervals. Occasionally there are exceptional components (often the one containing the trivial representation) whose homeomorphism types are not ``familiar'' in this sense, but our analysis nonetheless gives a stratification of each such component into more manageable pieces. The following three theorems, taken from Section \ref{sec:SomeMaps3}, are emblematic of the types of results appearing in Section \ref{sec:SpecialCases}.

\begin{theorem}[Surfaces]\label{thm:A}
Suppose $\Sigma$ has genus $g$ and $\pi = \bb{Z}_{2g+1}$ acts on $\Sigma$ with three fixed points. Assume two of these fixed points have the same weight. 
\begin{itemize}
\item The $\pi$-fixed point set on $\R(\Sigma)$ is $S^3$-equivariantly homeomorphic to a disjoint union of a point and $g$ copies of $S^2$, with $S^3$ acting on $S^2 \subseteq \bb{R}^3$ through the two-fold cover $S^3 \rightarrow \SO(3)$. 
\item The fixed point set on $\X(\Sigma)$ consists of $\frac{1}{6}(2g^3 - 3g^2 + 7g + 6)$ connected components. Of these, one is homeomorphic to an interval, $g$ are homeomorphic to $S^1$, and the remaining components are points.
\end{itemize}
\end{theorem}

\begin{theorem}[Circle bundles]\label{thm:C}
With $\Sigma$ and $\pi$ as in Theorem \ref{thm:A}, suppose $N \rightarrow \Sigma$ is a $\pi$-equivariant oriented circle bundle with Euler number $k$, and assume $\pi$ acts freely on $N$. The character variety $\X(N/ \pi)$ is homeomorphic to a disjoint union of an interval and $2g\lfloor \frac{k+1}{2} \rfloor$ points. 
\end{theorem}

\begin{theorem}[Mapping tori]\label{thm:B}
With $\Sigma$ and $\pi$ as in Theorem \ref{thm:A}, let $\phi \in \pi$ be a generator. The character variety $\X(\Sigma_\phi)$ has $\frac{1}{3}(2g+1)(g^2 - 2g + 3)$ connected components. Of these, one is homeomorphic to an interval, $g$ are homeomorphic to $S^1$, and the remaining components are points.
\end{theorem}

Other families of actions considered in Section \ref{sec:SpecialCases} include hyperelliptic involutions and certain rotations. While this section gives a deep dive into various special families, Section \ref{sec:InGeneral} gives a general approach to handling all such actions on surfaces, provided there is at least one fixed point. In this case, $\pi \cong \bb{Z}_n$ is necessarily cyclic. To simplify the discussion, we restrict to the case where the quotient $\Sigma / \bb{Z}_n$ has minimal genus, a condition equivalent to requiring that the genus of $\Sigma$ be $(t- 2)(n-1)/2$, where $t \geq 2$ is the number of fixed points; see Remark \ref{rem:genus} (a) for extending to higher genus quotients. Our starting point for this general discussion is Theorem \ref{thm:Action}, which gives a presentation for $\pi_1(\Sigma)$ on which the $\bb{Z}_n$-action is given by permuting the generators in a simple manner, and the relations are given in a $\bb{Z}_n$-invariant fashion. This may be of independent interest, so we reproduce the statement here:

\begin{theorem*}[Thm. \ref{thm:Action}]
With $\Sigma$ and $\pi \cong \bb{Z}_n$ as above, let $k_1, \ldots, k_t$ be the weights of the $\pi$-action on $\Sigma$. The fundamental group $\pi_1(\Sigma)$ has a presentation with generators $\alpha_{j, \ell}$ for $2 \leq j \leq t$ and $\ell \in \bb{Z}_n$, subject to the relations $\prod_{i = 0}^{n-1} \alpha_{j, i k_1} = e$ for $2 \leq j \leq t $ and
$$\prod_{j = 2}^t \prod_{i = 1}^{k^{*}_j} \alpha_{j, \ell(i, j)+ \mu} = e$$
for $\mu \in \bb{Z}_n$, with $\ell(i, j)$ as in (\ref{eq:elldef}) and $k_j^*$ the multiplicative inverse in $\bb{Z}_n$ of $k_j$. The action of $1 \in \bb{Z}_n$ on the generators sends $\alpha_{j, \ell} \mapsto \alpha_{j, \ell - k_1}$. 
\end{theorem*}

This result is predicated on Nielson's classification \cite{Nielson} for periodic actions, and we give a detailed proof of Theorem \ref{thm:Action} in Section \ref{sec:ProofOfTheoremAction}. Using Theorem \ref{thm:Action}, we give a description of the fixed point sets on $\R(\Sigma)$ and $\R(N)$; see Propositions \ref{prop:RepGen} and \ref{prop:RepNGen}. As a consequence we will see that, somewhat surprisingly, the topological type of each of these fixed point sets is independent of the choice of weights. For example, the description in Theorem \ref{thm:A} of the $\pi$-fixed point set in $\R(\Sigma)$ holds \emph{in the absence of} the stated hypothesis that at least two fixed points have the same weight. 

Despite giving a general description of the fixed point sets in $\R$, Propositions \ref{prop:RepGen} and \ref{prop:RepNGen} fall short of giving the homeomorphism types in as explicit a fashion as, say, our descriptions in Section \ref{sec:SpecialCases}. Nonetheless, in Section \ref{sec:UpTo4} we give such a description for $t  \leq 4$. We leave to future work an analogous description for $t \geq 5$, but in Section \ref{sec:t=5} we give a sense for the complications that arise for larger values of $t$ by computing the fixed point sets explicitly for $t = 5$ and $n = 3$.

\medskip

Finally, we point out a connection with Seifert-fibered spaces: Suppose $N \rightarrow \Sigma$ is a $\pi$-equivariant oriented circle bundle. When $\pi$ acts freely on $N$, the $\pi$-fixed points in $\R(N)$ can be identified with the elements of the representation variety $\R(N/ \pi)$ of the quotient (\emph{cf}. Proposition \ref{prop:SF}). The quotient $N/ \pi$ is naturally a Seifert-fibered space and so our computation of $\R(N/ \pi)$ complements the work of Fintushel--Stern \cite{FS} and Kirk--Klassen \cite{KK} who handle the case of Seifert-fibered \emph{homology spheres} (our spaces $\pi_1$-surject onto $\pi \cong \bb{Z}_n$ and so are not homology 3-spheres). We can already see in Theorem \ref{thm:C} a qualitative difference between our setting and the homology 3-sphere case: The character variety $\X(N/ \pi)$ has positive-dimensional components, while it is shown in \cite{FS} that the character variety of a Seifert-fibered homology 3-sphere is discrete. This discussion has extensions to the case where $\pi$ does not act freely on $N$, in which case the fixed point set on $\R(N / \pi)$ corresponds to the $S^3$-representations of the \emph{orbifold} fundamental group of $N/ \pi$; see Remark \ref{rem:Orbifolds}.

\medskip

\textbf{AI (non)-usage statement}: Generative AI was not used in any way in the proofs, results, writing, or
production of this paper. 

\medskip

\textbf{Acknowledgements}: This material is based upon work supported by the National Science Foundation (Grant No. 2349593) and the Elizabeth McConnell Bliss Endowment.

\section{Set-up and background}

Suppose $\pi$ is a finite group acting on a topological space $X$ through homeomorphisms. We assume all group actions are effective. We will write $\fix(\pi, X)$ for the fixed point set in $X$ of this action and $\fix(\phi, X)$ for the fixed point set of any $\phi \in \pi$. When $X$ is not directly relevant, or is clear from context, we drop it from the notation. Of course, $\fix(\pi) = \cap_{\phi \in \pi} \fix(\phi)$. If $\pi$ is cyclic and $\phi$ is a generator, then $\fix(\pi) = \fix(\phi)$. If $\pi$ is cyclic, then the $\pi$-action is semifree if and only if $\fix(\pi) = \fix(\phi)$ for all nonidentity elements $\phi \in \pi$; when $\pi$ has prime order the action is automatically semifree.

In all of our examples  $\pi \cong \bb{Z}_n$ will be cyclic, so we can fix a generator $\phi \in \pi$. When $X = \Sigma$ is a surface, from this choice of generator we can associate to each $\pi$-fixed point $x \in \Sigma$ a weight as follows: There is an orientation-preserving embedding $\bb{D}^2 \rightarrow \Sigma$, taking $0$ to $x$ so that, in these coordinates, the map $\phi$ takes the form $\bb{D}^2 \ni z \mapsto e^{2 \pi \hat{i} k / n} z$ for a unique $k \in \bb{Z}_n$; this $k$ is the \emph{weight} (or \emph{rotation number}) of $x$. Throughout we write $\hat{i} \defeq \sqrt{-1}$ to free up the symbol $i$ as an index. Since our group actions are effective and $\Sigma$ is connected, the weights will always be units in $\bb{Z}_n$, and all fixed points will be isolated.

\subsection{Quaternions and \texorpdfstring{$S^3$}{}}

We will view the space $\bb{H} = \bb{C} \oplus \bb{C} \hat{j}$ of quaternions as a vector space over $\bb{C}$ with basis $1, \hat{j}$. In this way, $\bb{C}$ is naturally a subspace of $\bb{H}$, as is $\bb{R}$, which is the center of $\bb{H}$. Every element of $\bb{H}$ can be conjugated (in the group-theoretic sense) to an element of $\bb{C}$, though this element is generally not unique: for each $z \in \bb{C}$ there is some $A \in \bb{H}$ so that $A z A^{-1} = \overline{z}$ is the complex conjugate of $z$. We will write $\mathrm{Re}(A) = \frac{1}{2}(A + A^*)$ and $\mathrm{Im}(A) = \frac{1}{2} (A - A^*)$ for the real and imaginary parts of $A$, respectively, where $A^*$ is the quaternionic conjugate of $A$. The non-zero elements of $\bb{H}$ act on $\bb{H}$ by (group-theoretic) conjugation, and the conjugacy classes are exactly the 3-planes given by the fibers of the map $\mathrm{Re}: \bb{H} \rightarrow \bb{R}$. We will work relative to the inner product
$$\langle A , B \rangle \defeq \mathrm{Re}(A B^*)$$
on $\bb{H}$; this is bi-invariant under the left- and right-multiplication actions of the unit sphere $S^3\subseteq \bb{H}$ on $\bb{H}$. 

The center of the unit sphere $S^3 \subseteq \bb{H}$ is $\bb{R} \cap S^3 = \left\{ \pm 1 \right\}$, and every element of $S^3$ can be conjugated to an element of $S^1 = \bb{C} \cap S^3$. If $A \in S^3$, then its conjugacy class is a 2-sphere in $S^3$ of radius $\sqrt{1 - \mathrm{Re}(A)^2}$. 

We will write 
$$\zeta_n = \left\{ A \in S^3 \; \big| \; A^n = 1 \right\}$$ 
for the set of $n$th roots of unity in $S^3$. Of course, if $A \in S^3$ is an $n$th root of unity, then all elements in the conjugacy class of $A$ are also $n$th roots of unity. We see from this that $\zeta_n$ has $\left\lfloor \frac{n}{2} \right\rfloor +1$ connected components of which one consists (resp. two consist) of a single point and the remaining are diffeomorphic to $S^2$ when $n$ is odd (resp. even). Indeed, by conjugation, we see that each component is represented uniquely by a \emph{complex} root of unity in the closed upper half plane.

We end this section with several lemmas that will be useful in our analysis below. For the first, for $A \in S^3$, consider the map $\mathcal{F}_A: S^3 \rightarrow S^3$ given by the commutator $\mathcal{F}_A(B) = [A , B] \defeq A B A^{-1} B^{-1} = AB A^* B^*$. 

\begin{lemma}\label{lem:centralizer}
The centralizer $C_A \in S^3$ of $A$ acts freely and transitively on the fibers of $\mathcal{F}_A$ by right multiplication.
\end{lemma}

\begin{proof}
For $C \in C_A$ we have
$$\mathcal{F}_A(BC) = A(BC)A^{-1}(BC)^{-1} = \mathcal{F}_A(B),$$
so we plainly have an action on the fibers, and it is immediate that this action is free. To see it is transitive on the fibers, note that if $\mathcal{F}_A(B_0) = \mathcal{F}_A(B_1)$, then $B_0^{-1} B_1 \in C_A$ is in the center of $A$. 
\end{proof}

Let $P_A$ be the 3-plane in $\bb{H}$ through 1 and normal to $A$. 

\begin{lemma}\label{lem:image}
For $A \in S^3$, the image of $\mathcal{F}_A$ is the intersection $P_A \cap S^3$. This intersection is a 2-sphere if and only if $A \neq 1$, otherwise it is a point. 
\end{lemma}

\begin{proof}
The lemma is straight-forward when $A = 1$, so assume $A \neq 1$. If $X = \mathcal{F}_A(B)$ lies in the image of $\mathcal{F}_A$, then $BA^{-1} B^{-1} = A^{-1} X$. Taking the real part and rearranging, we see that
$$\langle (1 - X), A \rangle = \mathrm{Re}((1-X)A^*) = 0,$$
so $X \in P_A$; taking the norm of $X = A B A^* B^*$ shows $X \in S^3$. 

That the intersection $P_A \cap S^3$ is a topological 2-sphere follows from the geometric description of the 3-plane. Thus, $\mathcal{F}_A$ is a map from $S^3$ to the 2-sphere $P_A \cap S^3$. By Lemma \ref{lem:centralizer}, this descends to an embedding
$$S^3 / C_A \hookrightarrow P_A \cap S^3.$$
The quotient $S^3 / C_A$ is a 2-sphere (this is the Hopf fibration), so we have an embedding of a 2-sphere into a 2-sphere. This embedding is surjective by general topological considerations, so the image of $\mathcal{F}_A$ is all of $P_A \cap S^3$. 
\end{proof}

Allowing the planes $P_A$ to vary for $A \in S^1$, we obtain the following.

\begin{cor}\label{cor:imageofF_A}
Each element of $S^3$ is in the image of $\mathcal{F}_A$ for some $A \in S^1$. 
\end{cor}

Our last result of this section deals with roots of unity.

\begin{lemma}\label{lem:LiorsCountingLemma}
Let $g \geq 0$ be an integer. The set $\mathcal{C}$ of pairs $(A, B) \in S^3 \times S^3$ with $A^{2g+1} = B^{2g+1} = (AB)^{2g+1} = 1$ has $\frac{1}{6} (g+1)(g+2)(2g+3)$ connected components. The action of $S^3$ on $\mathcal{C}$ by diagonal conjugation restricts to each component to be a transitive action. Exactly one of these components is a point and consists of $(1, 1)$; $2g(g+1)$ of these components are diffeomorphic to $S^2$, and each element is a commuting pair; each of the remaining $\frac{1}{6} g (g-1)(2g-1)$ components is diffeomorphic to $\SO(3)$ and consists entirely of non-commuting pairs. 
\end{lemma}

We note that this lemma can alternatively be articulated in the framework of linkages, which were used by Fintushel--Stern \cite{FS} and Kirk--Klassen \cite{KK} to study the topology of the character variety of Seifert-fibered homology spheres. 

\begin{proof}[Proof of Lemma \ref{lem:LiorsCountingLemma}]
Set $n = 2g+1$, so $\mathcal{C}$ is the set of pairs $(A, B) \in \zeta_n \times \zeta_n$ with $AB \in \zeta_n$. Since $S^3$ is connected, we automatically have that the specified action of $S^3$ on $\mathcal{C}$ restricts to an action on each component. Our proof below will show that this action is transitive on each component, making each component a homogenous space for $S^3$. As such, the topological type of each component will be immediate from any isotropy subgroup; for example, if $(A, B)$ do not commute, then the isotropy subgroup is $\left\{ \pm 1 \right\} \subseteq S^3$, so it will follow that these belong to a component equivariantly diffeomorphic to $S^3 /  \left\{ \pm 1 \right\} \cong \SO(3)$. We will see that the only isotropy groups that arise are $S^3, S^1$, and $\left\{ \pm 1 \right\}$.

To analyze the components further, consider the projection $(A, B) \mapsto A$, which restricts to a well-defined map of the form $\Pi: \mathcal{C} \rightarrow \zeta_n$. The fiber over $1$ is $\left\{1 \right\} \times \zeta_n$, which consists of a single point $(1, 1)$ and $g$ copies of $S^2$, on which $S^3$ acts transitively. All pairs in $\left\{1 \right\} \times \zeta_n$ plainly commute. It thus suffices to determine the components of $\mathcal{C} \backslash \Pi^{-1}(1)$. Towards this end, we will analyze $\Pi^{-1}(A)$ for $A  \in \zeta_n \backslash \left\{ 1\right\}$. We can simultaneously conjugate so that $A = e^{2 \pi \hat{i} m / n}$ for some $m$ in the range $ 1 \leq m \leq g= \lfloor n/ 2 \rfloor$. When we wish to emphasize the dependence on $m$, we will write $A_m \defeq A$. Note that the centralizer $C_A = S^1$ is the circle. To see that $S^3$ acts transitively on the components of $\mathcal{C}$, it thus suffices to show that $C_A$ acts transitively on the components of $\Pi^{-1}(A)$.

The map $(A, B) \mapsto AB$ induces an isometry
$$\Pi^{-1}(A) \cong   \zeta_n \cap (A \zeta_n)$$
with the intersection of $\zeta_n$ and $A \zeta_n = \left\{A B \; \vert \; B \in \zeta_n \right\}$; moreover, this isometry is $C_A$-equivariant relative to the conjugation action of $C_A$ on both sides. This shifts our focus to $\zeta_n \cap (A \zeta_n)$. Write $p: \zeta_n \rightarrow \bb{C}$ for the restriction to $\zeta_n$ of the projection $ P: \bb{H} = \bb{C} \oplus \bb{C} \hat{j} \rightarrow \bb{C}$. The image of $p$ is the intersection of the unit disk $\bb{D}^2$ with the vertical lines $\cup_{0 \leq k \leq g} \cos(2 \pi k / n) \times \hat{i} \bb{R}$. See Figure \ref{fig:1} (left) for an illustration of $\mathrm{im}(p) \subseteq \bb{D}^2$.

If $z \in p(\zeta_n)$ lies in this image and lies in the interior of $\bb{D}^2$, then the preimage $p^{-1}(z)$ is a circle and no points in this preimage commute with $A$; if $z \in p(\zeta_n)$ lies in the boundary $\partial \bb{D}^2$, then the preimage $p^{-1}(z) = \left\{z \right\}$ consists of a single point and this commutes with $A$ ($p$ restricts to the ``height function'' on each 2-sphere component of $\zeta_n$). The circle $S^1 = C_A$ acts on each fiber of $p$ transitively by conjugation.  

The set $A\zeta_n$ is obtained from $\zeta_n \subseteq \bb{C} \oplus \bb{C} \hat{j}$ by rotating by angle $2 \pi m / n$ jointly in the $\bb{C}$-plane and the $\bb{C} \hat{j}$-plane. Thus, writing $p_A$ for the restriction of $P$ to $A \zeta_n$, the image of $p_A$ is the intersection with $\bb{D}^2$ of the lines 
$$e^{2 \pi \hat{i} m / n}(\cos(2 \pi k / n) \times \hat{i} \bb{R})$$ 
for $0 \leq k \leq g $; see Figure \ref{fig:1} (middle). The fibers of $p_A$ have the same description as those of $p$, with the circle $S^1$ still acting transitively on them by conjugation. It follows that the connected components of the intersection $\zeta_n \cap (A \zeta_n) $ are in one-to-one correspondence with the components of $\mathrm{im}(p) \cap \mathrm{im}(p_A)  \in \bb{D}^2$; since $A \neq 1$, the components of this intersection are isolated points. The foregoing discussion shows that $S^3$ acts transitively on each fiber of $p$, so to finish the proof it suffices to count the number of these intersection points in $\bb{D}^2$. 

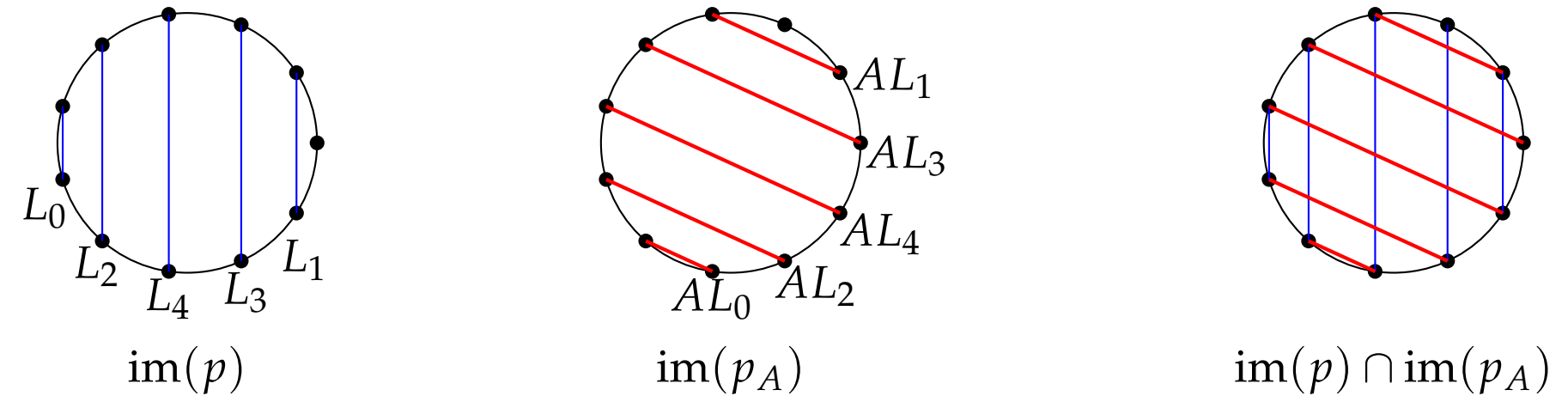
\begin{figure}[ht]
    \centering
\begin{tikzpicture}[
dot/.style={draw,fill,circle,inner sep=1pt}
]
\draw (0,0) circle (1);
\foreach \x in {0, ..., 10}, { 
 \node[dot] (\x) at ( \x*360/11:1) {} ;};
 \foreach \x in {1, ..., 5}, { 
\draw[-, blue] ( \x*360/11:1) -- ( -\x*360/11:1);
};
\node at (-1.1,-.5) {$L_0$};
\node at (-.7,-.95) {$L_2$};
\node at (-.15,-1.2) {$L_4$};
\node at (.45,-1.15) {$L_3$};
\node at (.9,-.9) {$L_1$};
        \node at (0,-1.75) {$\mathrm{im}(p)$};
  \end{tikzpicture}\hfill
   \begin{tikzpicture}[
dot/.style={draw,fill,circle,inner sep=1pt}
]
\draw (0,0) circle (1);
\foreach \x in {0, ..., 10}, { 
 \node[dot] (\x) at ( \x*360/11:1) {} ;};
 \foreach \x in {1, ..., 5}, { 
\draw[-, thick, red] (\x*360/11+720/11:1) -- ( -\x*360/11+720/11:1);};
         \node at (0,-1.75) {$\mathrm{im}(p_A)$};
\node at (-.15,-1.2) {$A L_0$};
\node at (.65,-1.1) {$A L_2$};
\node at (1.15,-.7) {$A L_4$};
\node at (1.35,-.1) {$A L_3$};
\node at (1.25,.5) {$A L_1$};
  \end{tikzpicture}\hfill
 \begin{tikzpicture}[
dot/.style={draw,fill,circle,inner sep=1pt}
]
\draw (0,0) circle (1);
\foreach \x in {0, ..., 10}, { 
 \node[dot] (\x) at ( \x*360/11:1) {} ;};
 \foreach \x in {1, ..., 5}, { 
\draw[-, blue] ( \x*360/11:1) -- ( -\x*360/11:1);};
 \foreach \x in {1, ..., 5}, { 
\draw[-, thick, red] (\x*360/11+720/11:1) -- ( -\x*360/11+720/11:1);};
         \node at (0,-1.75) {$\mathrm{im}(p) \cap \mathrm{im}(p_A)$};
  \end{tikzpicture}
  \hfill
  \caption{Illustrated here is the case $n = 11$ ($g = 5$) and $m = 2$. \hfill $\diamond$}
\label{fig:1}
\end{figure}

Of these intersection points, $n =2g+1$ lie on the boundary over which the fiber of $p$ consists of a single point that commutes with $A$. Recalling that $A$ is here parametrized by $\zeta_n \backslash \left\{1 \right\}$, which has $g$ connected components, the intersection points of this type produce $g(2g+1)$ components in $\mathcal{C}$ (each of which consists of commuting pairs). 

Moving now to the intersection points of $\mathrm{im}(p) \cap \mathrm{im}(p_{A_m}) $ lying in the interior of $\bb{D}^2$, note first that these all correspond to noncommuting pairs in $\mathcal{C}$. To finish the proof we need to count these points for $1 \leq m \leq g$. To carry this out, label the (positive-length) line segments of $\mathrm{im}(p)$ as $L_\ell$ for $ 0 \leq \ell \leq g-1$, as shown in Figure \ref{fig:1} (left). The defining feature of $L_\ell$ is that $\ell$ is the number of $n$th roots of unity contained in the interior of the smaller circular arc in $S^1$ defined by $L_\ell$. For $0 \leq \ell \leq g-1$, we will show that 
$$\sumdd{m = 1}{g} \big| \mathrm{im}(p) \cap A_m L_\ell \cap \mathrm{int}(\bb{D}^2) \big| = \frac{\ell(2g -\ell - 1)}{2};$$
summing over $\ell$ gives the number of components of this type claimed in the lemma. 

First assume $\ell$ is odd, so $L_\ell$ lies in the right-half plane. The line segment $A_1 L_\ell$ is obtained from $L_\ell$ by rotation by $2 \pi / n$. Then $A_1 L_\ell \cap \mathrm{im}(p)$ contains exactly one element in the interior of $\bb{D}^2$. Likewise $A_2 L_\ell \cap \mathrm{im}(p)$ contains $3$ and $A_3 L_\ell \cap \mathrm{im}(p)$ contains $5$, and so on. See Figure \ref{fig:2}. This pattern does not go on indefinitely: the maximum number of interior intersection points of $A_m L_\ell \cap  \mathrm{im}(p)$ is $\ell$, and this is realized for the first time when $m = (\ell + 1)/2$. This maximum continues to be realized as $m$ increases until $m= (\ell + 1)/2 + g-\ell - 1$. At the next value of $m$ (when $m =  (\ell + 1)/2 + g-\ell$), the number of interior intersection points of $A_m L_\ell \cap  \mathrm{im}(p)$ is $\ell -1$. At the next value this number is $\ell - 3$, then $\ell - 5$, a pattern that continues until $m = g$ when there are zero interior intersection points. In total, this gives 
$$1 + 2 + 3 + \ldots + (\ell-1) + (g- \ell) \ell = \frac{\ell(2g -\ell - 1)}{2}$$
interior intersection points, as claimed. A similar argument holds in the case where $\ell$ is even, the only significant difference being that the maximum is first realized when $m = \ell/2 + 1$ and last realized when $m = \ell / 2 + (g-\ell)$.

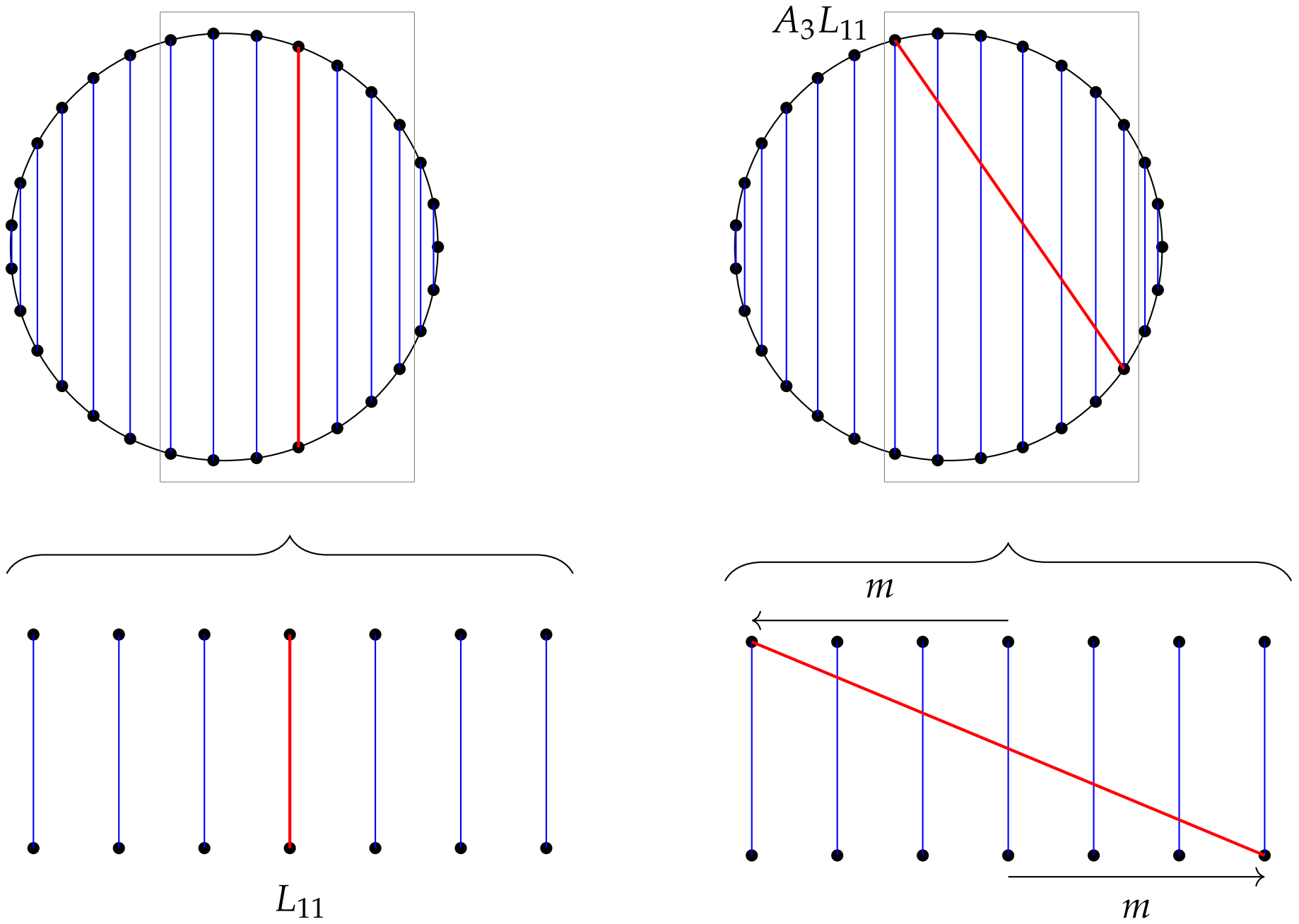
\begin{figure}[ht]
    \centering
\hspace{-1.25cm}
\begin{tikzpicture}[
dot/.style={draw,fill,circle,inner sep=1pt}
]
\draw (0,0) circle (2);
\draw[-, help lines] ( -.6, 2.2) -- ( 1.78, 2.2);
\draw[-, help lines] ( -.6, -2.2) -- ( 1.78, -2.2);
\draw[-, help lines] ( -.6, 2.2) -- ( -.6, -2.2);
\draw[-, help lines] ( 1.78, -2.2) -- ( 1.78, 2.2);
\foreach \x in {0, ..., 30}, { 
 \node[dot] (\x) at ( \x*360/31:2) {} ;};
 \foreach \x in {1, ..., 15}, { 
\draw[-, blue] ( \x*360/31:2) -- ( -\x*360/31:2);};
\draw[-, thick, red] ( 6*360/31:2) -- ( -6*360/31:2);
    \end{tikzpicture}  \hspace{2.5cm}
\begin{tikzpicture}[
dot/.style={draw,fill,circle,inner sep=1pt}
]
\draw (0,0) circle (2);
\draw[-, help lines] ( -.6, 2.2) -- ( 1.78, 2.2);
\draw[-, help lines] ( -.6, -2.2) -- ( 1.78, -2.2);
\draw[-, help lines] ( -.6, 2.2) -- ( -.6, -2.2);
\draw[-, help lines] ( 1.78, -2.2) -- ( 1.78, 2.2);
\foreach \x in {0, ..., 30}, { 
 \node[dot] (\x) at ( \x*360/31:2) {} ;};
 \foreach \x in {1, ..., 15}, { 
\draw[-, blue] ( \x*360/31:2) -- ( -\x*360/31:2);};
\draw[-, thick, red] ( 9*360/31:2) -- ( -3*360/31:2);
\node at (-1.2,2.1) {$A_{3}L_{11}$};
    \end{tikzpicture} \hspace{.5cm} \\
 \begin{tikzpicture}[
dot/.style={draw,fill,circle,inner sep=1pt}
]
\draw [ decorate,decoration={brace,amplitude=10pt},yshift=2pt] (-.25,2.5) -- (6*.8+.25,2.5) node [black, midway,yshift=20pt]{};
\foreach \x in {0, ..., 6}, { 
 \node[dot] (\x) at (\x*.8, 2) {} ;
 \node[dot] (\x) at (\x*.8, 0) {} ;};
 \foreach \x in {0, 1, 2, 4, 5, 6}, { 
\draw[-, blue] ( \x*.8,0) -- (\x*.8, 2);};
\draw[-, thick, red] ( 3*.8,0) -- (3*.8, 2);
\node at (2.5,-.5) {$L_{11}$};
  \end{tikzpicture}
  \hfill \hfill
 \begin{tikzpicture}[
dot/.style={draw,fill,circle,inner sep=1pt}
]
\draw [ decorate,decoration={brace,amplitude=10pt},yshift=2pt] (-.25,2.5) -- (6*.8+.25,2.5) node [black,midway,yshift=20pt]{};
\foreach \x in {0, ..., 6}, { 
 \node[dot] (\x) at (\x*.8, 2) {} ;
 \node[dot] (\x) at (\x*.8, 0) {} ;};
 \foreach \x in {0, ..., 6}, { 
\draw[-, blue] ( \x*.8,0) -- (\x*.8, 2);};
\draw[-, thick, red] ( 6*.8,0) -- (0*.8, 2);
\draw[->] ( 3*.8,-.2) -- (6*.8, -.2);
\draw[->] ( 3*.8,2.2) -- (0*.8, 2.2);
\node at (.8*3*6/4,-.5) {$m$};
\node at (.8*6/4,2.5) {$m$};
  \end{tikzpicture}
  \hfill
    \caption{This illustrates the case $g = 15$, $m = 3$, and $\ell = 11$. On the right, the angled (thick, red) line now intersects $ 5  = 2m-1$ segments in the interior of the disk. As $m$ increases, the number of intersection points increases by 2, until the maximum number of 11 intersection points is reached. Once the left endpoint of $A_m L_{11}$ drops below the $x$-axis, the number of intersection points for $A_m L_{11}$ decreases to zero at $m = g$, since $A_{g} L_{11}$ is nearly vertical. \hfill $\diamond$}
    \label{fig:2}
\end{figure}
\end{proof}

\subsection{The representation and character varieties}

We endow $\R(X) = \hom(\pi_1(X, x_0), S^3)$ with the compact-open topology, relative to the discrete topology on $\pi_1(X, x_0)$. When the basepoint $x_0$ is relevant to the discussion, we will write $\R(X, x_0)$. Our spaces $X$ will be path-connected, so $\R(X, x_0)$ is independent of $x_0$ up to homeomorphism. 

Our spaces $X$ will also be compact, so $\pi_1(X, x_0)$ will be finitely-presented. Given a presentation 
$$\pi_1(X, x_0) = \langle \gamma_1, \ldots, \gamma_k \vert w_1, \ldots, w_\ell \rangle$$
the evaluation map
\begin{equation}\label{eq:Rhomeo}
\R(X) \longrightarrow (S^3)^k, \hspace{1cm} \rho \longmapsto (\rho(\gamma_1), \ldots, \rho(\gamma_k))
\end{equation}
is injective with image given by the tuples $(g_1, \ldots, g_k) \in (S^3)^k$ satisfying $W_1 = \ldots = W_\ell = 1$, where $W_j$ is exactly the word $w_j$, but with $g_i$ in place of $\gamma_i$. Importantly, with the specified topology on the representation variety, the map (\ref{eq:Rhomeo}) is a homeomorphism onto its image, and so we are free to identify $R(X)$ with its image, provided the presentation is clear from context. Moreover, $\R(X)$ is locally path-connected, so the path components are precisely the connected components. 

For example, when considering a surface $\Sigma = \Sigma_g$ of genus $g \geq 1$, we may choose to work with a generating set $\alpha_1, \beta_1, \ldots, \alpha_g, \beta_g$ for $\pi_1(\Sigma, x_0)$, that is \emph{standard} in the sense that 
$$\prod_{j = 1}^g [\alpha_j, \beta_j ] =e$$
where $[\alpha, \beta] = \alpha \beta \alpha^{-1} \beta^{-1}$ is the group commutator and $e$ is the identity. Here and below the product notation refers to the increasing ordering convention: if $b_1, \ldots b_n$ are elements of a group, then
$$\prod_{i = 1}^{n} b_i \defeq b_1 b_2 b_3 \ldots b_n.$$ 
In this case, the map (\ref{eq:Rhomeo}) gives a homeomorphism
$$\R(\Sigma) \cong \left\{(A_1, B_1, \ldots, A_g, B_g) \in (S^3)^{2g}\; \Bigg| \; \prod_{j = 1}^g [A_g, B_g] = 1\right\}. $$
The representation variety of $\Sigma$ is connected, which is an observation that can be shown without much trouble starting with the observation that the commutator map $S^3 \times S^3 \rightarrow S^3$ is surjective with connected fibers. 

The group $S^3$ acts on $\R(X)$ by conjugation on the codomain of our representations. Since $S^3$ is connected, the number of connected components of the character variety $\X(X) = \R(X)/S^3$ equals that of $\R(X)$:
\begin{equation}\label{eq:pi0}
\pi_0(\R(X)) \cong \pi_0(\X(X)).
\end{equation}

A representation $\rho \in \R(X)$ is \emph{reducible} if its image is contained in a maximal torus of $S^3$, and \emph{central} if its image is contained in the center of $S^3$. The reducible and central representations often are singular points in the representation variety. The \emph{trivial} representation is the one sending all elements of $\pi_1(X)$ to $1 \in S^3$.

If $\phi: X \rightarrow Y$ is continuous, then precomposition by $\phi_*: \pi_1(X) \rightarrow \pi_1(Y)$ determines a \emph{pullback} map
$$\phi^*: \R(Y, \phi(x_0)) \longrightarrow \R(X, x_0), \hspace{.5cm} \rho \longmapsto \rho \circ \phi_*$$
and this is continuous. The pullback descends to a continuous map
$$\phi^*: \X(Y) \longrightarrow \X(X)$$
on the character variety that we will denote by the same symbol. 

When $X = Y$ and $x_0$ is a fixed point of $\phi$, we introduce the following shorthand:
$$\fix_\R(\phi) \defeq \fix(\phi^*, \R(X)), \hspace{.5cm} \fix_\X(\phi) \defeq \fix(\phi^*, \X(X)).$$
More generally, if $\pi \leq \homeo(X)$ is a subgroup that fixes our basepoint $x_0$, then pullback by elements of $\pi$ induces an action of $\pi$ on $\R(X)$ and on $\X(X)$ that we still denote by $\pi$, and we set
$$\fix_\R(\pi) \defeq \fix(\pi, \R(X)), \hspace{.5cm} \fix_\X(\pi) \defeq \fix(\pi, \X(X)).$$

\subsection{Circle bundles and finite group actions}

Let $\Sigma = \Sigma_g$ be a genus-$g$ surface and suppose $N \rightarrow \Sigma$ is an oriented circle bundle with Euler number $k \in \bb{Z}$. When we wish to express the genus or Euler class explicitly, we may write $N_g, N^k$, or $N_g^k$ in place of $N$. The fundamental group of $N$ has presentation with generators $\alpha_1, \beta_1, \ldots, \alpha_g, \beta_g, z$ with $z$ central and 
$$z^k = \prod_{j = 1}^g [\alpha_j, \beta_j ].$$
As such, the representation variety $\R(N)$ can be identified with the set of tuples $(A_1, B_1, \ldots, A_g, B_g, Z) \in (S^3)^{2g+1}$ satisfying
\begin{equation}\label{eq:pi1circlebundle}
Z^k = \prod_{j = 1}^g [A_j, B_j ], \hspace{.5cm} [A_j, Z] = [B_j, Z] = 1.
\end{equation}

Now assume $\pi$ is a cyclic group acting on $\Sigma$ with fixed point $x_0 \in \Sigma$. Assume further that the action of $\pi$ on $\Sigma$ lifts to an action of $\pi$ on $N$ making $N$ a $\pi$-equivariant bundle. When we view $\pi$ as acting on $N$, we are referring to a homomorphism $\pi \rightarrow \homeo(N)$. We will write $\widetilde{\phi}$ for the image of $\phi$ under this homomorphism, so $\widetilde{\phi}: N \rightarrow N$ is a homeomorphism covering the homeomorphism $\phi: \Sigma \rightarrow \Sigma$.

In general, this lifted action of $\pi$ on $N$ need not have a fixed point, and this complicates matters at the fundamental group level since the elements of $\pi$ need not take $\pi_1(N, y_0)$ to itself for any $y_0 \in N$. Nonetheless, we can rectify this as follows: First, take as a basepoint any point $y_0 \in N$ in the fiber over $x_0 \in \Sigma$. Though $\pi$ need not fix $y_0$, it will take $y_0$ to another element in this same fiber. Since the fiber over $x_0$ is an oriented circle, for each $\phi \in \pi$ there is a canonical path $\tau$ from $y_0$ to $\widetilde{\phi}(y_0)$ given by the positively oriented arc. This $\tau$ induces a change of basepoint isomorphism
\begin{equation}\label{eq:basepointiso}
\pi_1(N, \widetilde{\phi}(y_0)) \cong \pi_1(N, y_0).
\end{equation}
The homotopy class $z$ from the opening paragraph of this section is represented by the loop traversing this same fiber (in the positive direction) that starts and ends at $y_0$. Since $z$ is central in $\pi_1(N, y_0)$, we see that the isomorphism (\ref{eq:basepointiso}) is independent of our choice of path in the fiber from $y_0$ to $\widetilde{\phi}(y_0)$ (e.g., we could replace $\tau$ with any loop of the form $\tau * z^k$ and obtain the same isomorphism). As such, by appealing to the isomorphism (\ref{eq:basepointiso}), we obtain a well-defined action of $\pi$ on $\pi_1(N, y_0)$ that covers the action on $\pi_1(\Sigma, x_0)$ relative to the map $\pi_1(N, y_0) \rightarrow \pi_1(\Sigma, x_0)$ induced by the projection $N \rightarrow \Sigma$. Using brackets $[\cdot ]$ for based homotopy class and $*$ for path concatenation, this action at $[\eta] \in \pi_1(N, y_0)$ is given by
$$\widetilde{\phi} \cdot [\eta] = [\tau * (\widetilde{\phi} \circ \eta) * \tau^{-1} ].$$
 In particular, we get an action of $\pi$ on the representation variety $\R(N) = \hom(\pi_1(N, y_0), S^3)$.

\begin{remark}\label{rem:NRem}
(a) Any $\pi$-action on $N$ will necessarily preserve the center of $\pi_1(N, y_0)$; so any given element of $\pi$ maps $z$ to $z$ or $z^{-1}$. It is shown in \cite[Cor. 3.1]{CT} that, provided $(g, k) \neq (1, 0)$, the action of $\pi$ on $\pi_1(\Sigma)$ always lifts to an action of $\pi$ on $\pi_1(N)$ acting trivially on the center. To refer to lifts of this type, we will say this lift of the $\pi$-action to $\pi_1(N)$ \emph{preserves the fiber-orientation}. We note, however, that there is a technical topological problem we do not address, which is whether this lifted action of $\pi$ on $\pi_1(N)$ is induced from a $\pi$-action on $N$ through homeomorphisms (the ``Nielson realization problem'' for circle bundles). Instead we will simply assume that there is an action of $\pi$ on $N$ covering the action on $\Sigma$ and so that the induced action on $\pi_1(N)$ preserves the fiber-orientation; we will abbreviate this by saying that $N$ is a $\pi$-\emph{equivariant oriented circle bundle} over $\Sigma$.

\medskip

(b) Here is an extension of this construction: Let $\widetilde{N}$ be the universal cover of $N$, viewed as the based-homotopy classes $[\gamma]$ of paths $\gamma: [0, 1] \rightarrow N$ with $\gamma(0) = y_0$. For $\widetilde{\phi} \in \pi$, the composition $\widetilde{\phi} \circ \gamma$ is a path in $N$ starting at $\widetilde{\phi}(y_0)$. As above, we can concatenate this with the path $\tau$, and thus obtain a path in $N$ starting at $y_0$, which we will denote by $\widetilde{\phi} \cdot \gamma \defeq \tau * (\widetilde{\phi} \circ \gamma)$. Then the map $[\gamma] \mapsto \widetilde{\phi} \cdot [\gamma] \defeq [\widetilde{\phi} \cdot \gamma]$ induces an action of $\pi$ on $\widetilde{N}$. 

Relative to the action of $\pi_1(N, y_0)$ on $\widetilde{N}$ by path-concatenation, this satisfies
\begin{equation}\label{eq:fixedaction}
\widetilde{\phi}([\eta] * [\gamma]) = (\widetilde{\phi} \cdot [\eta]) * (\widetilde{\phi} \cdot [\gamma])
\end{equation}
for $[\eta] \in \pi_1(N, y_0)$ and $[\gamma] \in \widetilde{N}$. \hfill $\diamond$
\end{remark}

Though there may be multiple distinct $\pi$-actions on $N$ covering the same $\pi$-action on $\Sigma$, these different actions are not detected by the homeomorphism type of the representation or character varieties:

\begin{prop}
Suppose $N$ and $N'$ are two $\pi$-equivariant oriented circle bundles over $\Sigma$ covering the same $\pi$-action on $\Sigma$. Writing $k$ and $k'$ for the Euler numbers of $N$ and $N'$, respectively, if $\vert k \vert = \vert k' \vert$, then the fixed point sets
$$\fix(\pi, \R(N)) \cong \fix(\pi, \R(N')), \hspace{.5cm} \fix(\pi, \X(N)) \cong \fix(\pi, \X(N'))$$
are homeomorphic. 
\end{prop}

\begin{proof}
Let $\phi \in \pi$ be a generator. Using the description (\ref{eq:pi1circlebundle}) of the representation variety $\R(N)$, and the discussion above about the action of $\pi$ on the fundamental group of $N$, we have
$$\widetilde{\phi}^*(A_1, B_1, \ldots, A_g, Z) = (\phi^*(A_1, B_1, \ldots, A_g, B_g), Z).$$
The right-side is independent of this lift $\widetilde{\phi}$. Since (\ref{eq:pi1circlebundle}) is independent of the sign of $k$, the proposition follows.
\end{proof}

When $\pi$ acts freely on $N$, the quotient $N / \pi$ is a smooth 3-manifold fibered by circles, and so is a Seifert-fibered space. The following, whose proof we defer to Section \ref{sec:GaugeTheory}, is how the computations of the present paper can be used to compute the representation and character varieties of Seifert fibered spaces of this type. 

\begin{prop}\label{prop:SF}
When $\pi = \bb{Z}_n$ acts freely on $N$, the pullback by the projection $N \rightarrow N/ \pi$ induces homeomorphisms
$$\R(N/ \pi) \cong \fix( \pi, \R(N)), \hspace{.5cm} \X(N/ \pi) \cong \fix( \pi, \R(N))/ S^3.$$
\end{prop}

See Remark \ref{rem:Orbifolds} for an extension without the requirement that the $\pi$-action is free.

\subsection{Mapping tori}

The \emph{mapping torus} of a homeomorphipsms $\phi: X \rightarrow X$ is the space
\[X_\phi = \frac{[0,1] \times X}{(1, x) \sim (0, \phi(x))}.\]
In our examples $\phi$ will always have a fixed point $x_0$, which we will take as the basepoint in $X$. In the mapping torus we will take as basepoint the equivalence class $\widetilde{x}_0$ of $(0, x_0)$. The fundamental group $\pi_1(X_\phi, \widetilde{x}_0)$ has a presentation
\[\langle \gamma, \tau \mid \phi_*\gamma = \tau^{-1} \gamma \tau, \;\; \forall \gamma \in \pi_1(X, x_0) \rangle \]
where $\tau$ can be viewed as the homotopy class of the image in $X_\phi$ of $t \mapsto (t, x_0)$. 

There is a close relationship between the character variety $\X(X_\phi)$ of the mapping torus and the fixed point set of $\phi^*$ on $\X(X)$. Indeed, pullback by the fiber inclusion $X \hookrightarrow X_\phi$ gives a surjection
\begin{equation}\label{eq:chimap}
\X(X_\phi) \longrightarrow \fix_\X(\phi)
\end{equation}
and the fiber over $[\rho]$ is can be identified with the centralizer $C_\rho \subseteq S^3$ of (the image of) any representative $\rho \in [\rho]$.

\subsection{Gauge theory}\label{sec:GaugeTheory}

Let $X$ be a smooth manifold; to simplify the discussion, we will assume $X$ is connected and of dimension at most 3. Write $\A(X)$ and $\G(X)$ for the spaces of connections and gauge transformations, respectively, on the trivial $S^3$-bundle over $X$ (\emph{cf}. \cite{KN}). Recall that a connection $A \in \A(X)$ is \emph{flat} if its curvature $F_A$ is identically zero; we will write $\A_{\flat}(X)$ for the subspace of flat connections. The gauge group acts on $\A(X)$ and this action preserves $\A_{\flat}(X)$. Let $x_0 \in X$ be the basepoint at which the fundamental group is defined and write $\G_0(X)$ for the subgroup consisting of $u \in \G(X)$ with $u(x_0) = 1$; we call $\G_0(X)$ the group of \emph{based} gauge transformations. The holonomy induces a map
$$H_X: \A_{\flat}(X) \longrightarrow \R(X)$$
that is invariant under the action of $\G_0(X)$. Moreover, this map is surjective and induces a homeomorphism
$$\frac{\A_{\flat}(X)}{\G_0(X)} \cong \R(X).$$
We briefly recall the surjectivity of $H_X$: Any $\rho \in \R(X)$ determines an action of $\pi_1(X, x_0)$ on $\widetilde{X} \times S^3$, where $\widetilde{X} \rightarrow X$ is the universal cover. The trivial connection on $\widetilde{X} \times S^3$ is $\pi_1(X, x_0)$-invariant and so descends to a well-defined flat connection $A_\rho$ on $P_\rho \defeq (\widetilde{X} \times S^3) / \pi_1(X, x_0)$. This quotient $P_\rho$ is a principal $S^3$-bundle over $X$ and so is necessarily trivializable since $\dim(X) \leq 3$. Upon choosing a trivialization, we can view $A_\rho$ as an element of $\A_{\flat}(X)$ and it follows that $H_X(A_\rho) = \rho$. 

We plainly have $\G(X) / \G_0(X) \cong S^3$ and so a consequence of the above construction is that $H_X$ induces a homeomorphism
$$\frac{\A_{\flat}(X)}{\G(X)} \cong \X(X)$$
with the character variety. 

\begin{proof}[Proof of Prop. \ref{prop:SF}]
In the notation of this proposition, write $\widetilde{\phi}: N \rightarrow N$ for a generator of $\pi$. Let $y_0 \in N$ be a point in the fiber over a fixed point on $\Sigma$, and write $[y_0]$ for the image of $y_0$ in $N/ \pi$. For $\rho \in \R(N)$, write $A_\rho$ for the flat connection on $P_\rho \rightarrow N$ as above. If $\rho \in \fix(\pi, \R(N))$ is fixed under the action of $\pi$, then it follows from (\ref{eq:fixedaction}) that the action of $\pi$ on $\widetilde{N}$ induces a natural action of $\pi$ on $P_\rho$ making $P_\rho \rightarrow N$ a $\pi$-equivariant principal $S^3$-bundle. Moreover, $A_\rho$ is $\pi$-invariant and so descends to a well-defined flat connection $A_\rho / \pi$ on $P_\rho / \pi \rightarrow N/ \pi$; since we have assumed $\pi$ acts freely on $N$, these quotients are smooth manifolds. In this way, we obtain a map
$$\fix(\pi, \R(N)) \longrightarrow \mathcal{A}_{\flat}(N/\pi) / \mathcal{G}_0(N/ \pi) \cong \R(N/ \pi).$$

For the reverse direction, use the fact that each element of $\R(N/ \pi)$ can be represented by a flat connection on $N / \pi$. By pulling this back to $N$, we get a $\pi$-invariant flat connection $A$ on $N$. This $A$ descends to an element of $\fix(\pi, \R(N))$ that is independent of the choices involved (since any two such $A$ differ by a $\pi$-invariant element of $\mathcal{G}_0(N)$). 

The claim about the character varieties is similar, but one should replace the based gauge group $\mathcal{G}_0(N)$ with the full gauge group $\mathcal{G}(N)$. 
\end{proof}

\begin{remark}\label{rem:Orbifolds}
There are extensions of Proposition \ref{prop:SF} to the case where $\pi$ does not act freely on $N$, but such extensions naturally involve orbifolds (V-manifolds). In this more general case, we still have that the elements of $\fix(\pi, \R(N))$ correspond to $\pi$-invariant flat connections on $N$, but now they descend to flat \emph{orbifold} connections on the orbifold $N/ \pi$. When taken modulo based gauge transformations, these are identified with the representations $\pi_1^{orb}(N / \pi) \rightarrow S^3$ of the orbifold fundamental group of $N / \pi$ \cite{FS}, but see also \cite[Thm. 3.4]{FuS}. In this way, we have a homeomorphism
$$\R^{orb}(N/ \pi)  \cong \fix(\pi, \R(N))$$
where $\R^{orb}(N/ \pi)$ is the set of homomorphisms $\pi_1^{orb}(N/ \pi) \rightarrow S^3$. A similar statement holds for the character variety. This gives an interpretation of the fixed point sets $\fix(\pi, \R(N))$ and $\fix(\pi, \X(N))$ even when $\pi$ need not act freely on $N$. \hfill $\diamond$
\end{remark}

\section{Significant special cases}\label{sec:SpecialCases}

Here we highlight several special cases and thoroughly analyze the relevant fixed point sets in each case.

\subsection{Hyperelliptic involutions (\texorpdfstring{$\pi = \bb{Z}_2$}{})}\label{sec:HEI}

Let $\Sigma$ be a surface of genus $g \geq 1$ and suppose $\phi: \Sigma \rightarrow \Sigma$ is an order-two map that is orientation-preserving and acts by multiplication by $-1$ on $H_1(\Sigma)$. It follows that $\phi_* \gamma = \gamma^{-1}$ for all $\gamma \in \pi_1(\Sigma)$. Write $\pi = \langle \phi \rangle \cong \bb{Z}_2$ for the subgroup of $\homeo(\Sigma)$ generated by $\phi$. 

\begin{figure}[h]
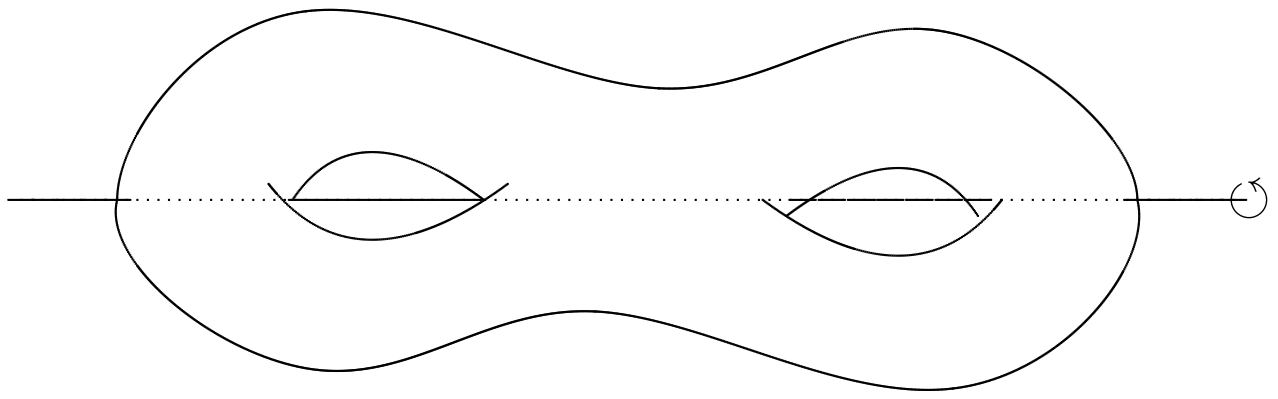
\label{figure0}
\centerline{\xy 
		/r.5pc/:,
(16,0)*{}="A";
(-16, 0)*{}="B";
"A";"B" **\crv{(16,2)& (9,8)& (2,1)& ( -11,9)& (-16,2)};  
"A";"B" **\crv{ (16.5, -2)& (11,-9)& (-2,-1)& (-9,-8)& (-16.5,-2)};  
(-4.5,0)*{}="b1";
(-10.5,0)*{}="b2";
"b1";"b2" **\crv{(-8.5,3)}; 
(-3.75,.5)*{}="B1";
(-11.25,.5)*{}="B2";
"B1";"B2" **\crv{ (-8.5,-3)}; 
(5,-0.5)*{}="t1";
(11,-0.5)*{}="t2";
"t1";"t2" **\crv{(9,2.5)}; 
(4.25,0)*{}="T1";
(11.75,0)*{}="T2";
"T1";"T2" **\crv{ (9,-3.5)}; 
(19,0)*{}="A1";
"A1"; "A" **\crv{};
(-19,0)*{}="B1";
"B1"; "B" **\crv{};
(-10.2,0)*{}="C1";
(-5,0)*{} = "C";
"C1"; "C" **\crv{};
(11,0)*{}="D1";
(5.5,0)*{} = "D";
"D1"; "D" **\crv{};
"A1"; "B1" **\crv{~**\dir{.}};
(19.5,0)*{\circlearrowleft};
\endxy}
\caption{A hyperelliptic involution on a surface of genus 2. The involution spins the surface about the axis as indicated. \hfill $\diamond$}
\end{figure}

The fixed point set $\fix_\R(\bb{Z}_2)$ can be identified with the set of tuples 
$$(A_1, B_1, \ldots, A_g, B_g) \in (S^3)^{2g}$$
satisfying $A_j = A_j^{-1}, B_j = B_j^{-1}$ and $\prod_j [A_j, B_j] = 1$. Since the only order-two elements of $S^3$ are $\pm 1$, we have the following.

\begin{prop}
For $g \geq 0$, the fixed point set in $\R(\Sigma)$
$$\fix_\R(\bb{Z}_2) \cong \left\{ -1, 1 \right\}^{2g}$$
consists of $4^g$ points. 
\end{prop}

Now assume that $N = N^k \rightarrow \Sigma$ is a $\bb{Z}_2$-equivariant oriented circle bundle. Using the description (\ref{eq:pi1circlebundle}), we again have that $A_j, B_j \in \left\{ \pm 1 \right\}$ and we see that this forces $Z^k = 1$.

\begin{prop}
For $g \geq 0$, and $k \in \bb{Z}$, there is an $S^3$-equivariant homeomorphism
$$\fix(\bb{Z}_2, \R(N)) \cong \left\{ -1, 1 \right\}^{2g} \times \zeta_k$$
with $S^3$ acting on the right by diagonal conjugation.  
\end{prop}

Since $S^3$ acts transitively on each component of $\zeta_k$, we can combine this with Proposition \ref{prop:SF} to get the following. 

\begin{cor}
If the $\bb{Z}_2$-action on $N$ is free, then $\X(N/ \pi) $ consists of $4^g  (\lfloor k/2 \rfloor + 1)$ points.
\end{cor}

Next we wish the analyze the spaces $\R(\Sigma_\phi)$, $\X(\Sigma_\phi)$, and $\fix_\X(\phi)$. As will become apparent in our analysis below, the exact topological-type of these spaces is fairly complicated, so we will content ourselves with showing the following.

\begin{prop}\label{prop:HEIcon}
The spaces $\R(\Sigma_\phi)$, $\X(\Sigma_\phi)$, and $\fix_\X(\phi)$ are each connected. 
\end{prop}

\begin{proof}
In light of the surjections $\R(\Sigma_\phi) \rightarrow \X(\Sigma_\phi) \rightarrow \fix_\X(\phi)$, it suffices to show that $\R(\Sigma_\phi) $ is connected. This space consists of the tuples $((A_j, B_j)_{j= 1}^{2g}, T)$ satisfying $\prod_j [A_j, B_j] = 1$ and 
$$T^{-1} A_j T = A_{j}^{-1}, \hspace{.5cm}  T^{-1} B_j T = B_{j}^{-1}.$$
When any element of the $\left\{A_j, B_j \right\}_j$ is noncentral, then $T$ is necessarily purely imaginary and orthogonal to all elements of $\left\{A_j, B_j \right\}_j$. For $T \in S^3$ purely imaginary, we therefore consider the set
$$\mathcal{C}_T \defeq \left\{(A_j, B_j)_j\in (T^\perp )^{2g} \; \Bigg| \;   \prod_{j = 1}^g [A_j, B_j] = 1 \right\}.$$
We therefore see that
$$\R(\Sigma_\phi) = \Big( \left\{ \pm 1 \right\}^{4g} \times S^3 \Big) \cup \bigcup_{T \in S^3: \mathrm{Re}(T)= 0} \mathcal{C}_T \times \left\{ T \right\}.$$
Each element of $\left\{ \pm 1 \right\}^{4g} \times S^3$ can be path-connected through this space to $\mathcal{C}_{\hat{i}} \times \left\{ \hat{i}\right\}$, so the proposition is a consequence of the following claim:

\medskip

\noindent \emph{Claim: Each $\mathcal{C}_T$ is connected.}

\medskip

The case with $g = 1$ is straight-forward, since $\mathcal{C}_T$ consists of the commuting pairs in the 2-sphere $T^\perp \cap S^3$. We may therefore assume $g \geq 2$ and, by simultaneously conjugating, we may assume $T = \hat{j}$. For $X \in S^3$ define
$$\mathcal{C}_{T, X}(g) = \left\{ (A_j, B_j)_j \in (T^\perp \cap S^3)^{2g} \; \Bigg| \; \prod_{j = 1}^g [A_j, B_j] = X \right\},$$
(so $\mathcal{C}_T = \mathcal{C}_{T, 1}(g)$). We will prove by induction on $g \geq 2$ that (i) $\mathcal{C}_{T, 1}(g)$ is connected, (ii) $\mathcal{C}_{T, X}(g-1)$ is nonempty for each $X$ and (iii) $\mathcal{C}_{T, X}(g-1)$ varies continuously in $X$; the claim follows from this. 

To prove the base case, let $X \in S^3$. By the case with $g = 1$ just considered, we may assume $X \neq 1$. By Corollary \ref{cor:imageofF_A} there is some $A \in S^1 \subseteq T^\perp \cap S^3$ and some $B_0 \in S^3$ so that $[A, B_0] = X$. By Lemma \ref{lem:centralizer}, for each $C \in C_A$ we have $[A, B_0 C] = X$. The set $B_0 C_A$ is a great circle in $S^3$, and so necessarily intersects the (great) 2-sphere $T^\perp \cap S^3$; any such $B$ in this intersection gives a pair $(A, B) \in \mathcal{C}_{T, X}(1)$. In fact, since $X \neq 1$, this intersection consists of exactly two points that vary continuously with $X$. To see that $\mathcal{C}_{T, 1}(2)$ is connected, notice that $\mathcal{C}_{T, 1}(2)$ is the set of tuples $(A_1, B_1, A_2, B_2) $ so that 
$$(A_1, B_1) \in \mathcal{C}_{T, X}(1), \hspace{.5cm} \textrm{and} \hspace{.5cm} (A_2, B_2) \in \mathcal{C}_{T, X^{-1}}(1) $$
for some $X \in S^3$. Any element of this space can be path-connected to a tuple with $(A_1, B_1) \in \mathcal{C}_{T, 1}(1)$ and $(A_2, B_2) \in \mathcal{C}_{T, 1}(1)$. Since $\mathcal{C}_{T, 1}(1)$ is connected, the case for $g = 2$ follows. The inductive step is similar. 
\end{proof}

\subsection{A family of maps with three fixed points (\texorpdfstring{$\pi = \bb{Z}_{2g+1}$}{})}\label{sec:SomeMaps3}

\begin{figure}[ht]
    \centering
\begin{tikzpicture}[
dot/.style={draw,fill,circle,inner sep=1pt}
]
\begin{scope}[very thick,decoration={
    markings,
    mark=at position 0.5 with {\arrow{>}}}
    ] 
\foreach \x in {0, ..., 6}, {
\draw[postaction={decorate}] ( \x*360/7 - 360/7 + 180:2)--( \x*360/7 + 180:2);};
\end{scope}
\node (A1) at (0, 0) {\Huge $\circlearrowleft$};
\node at (1.7, 1.7) {\large $\alpha_1$};
\node at (-.5, 2.1) {\large $\alpha_2$};
\node at (-2.1, .8) {\large $\alpha_3$};
\node at (-2.1, -.8) {\large $\alpha_4$};
\node at (-.5, -2.1) {\large $\alpha_5$};
\node at (1.3,-1.8) {\large $\alpha_6$};
\node at (2.2, 0) {\large $\alpha_7$};
  \end{tikzpicture}
  \hspace{.5cm}
  \begin{tikzpicture}[
dot/.style={draw,fill,circle,inner sep=1pt}
]
\begin{scope}[very thick,decoration={
    markings,
    mark=at position 0.5 with {\arrow{>}}}
    ] 
\foreach \x in {0, ..., 6}, {
\draw[postaction={decorate}] (\x*360/7 + 360/7 :2) -- ( \x*360/7 :2);};
\end{scope}
\node (A1) at (0, 0) {\Huge $\circlearrowleft$};
\node at (-1.3,-1.8) {\large $\alpha_1$};
\node at (.5, -2.1) {\large $\alpha_2$};
\node at (2.1, -.8) {\large $\alpha_3$};
\node at (2.1, .8) {\large $\alpha_4$};
\node at (.5, 2.1) {\large $\alpha_5$};
\node at (-1.4,1.8) {\large $\alpha_6$};
\node at (-2.2, 0) {\large $\alpha_7$};
  \end{tikzpicture}
  \caption{Indicated here is the case $g = 3$. The two 7-gons glue together as indicated to produce $\Sigma_3$, and the rotation by $2 \pi / 7$ on each $7$-gon induces an action of $\bb{Z}_7$ on $\Sigma_3$.\hfill $\diamond$}
\label{fig:4}
\end{figure}
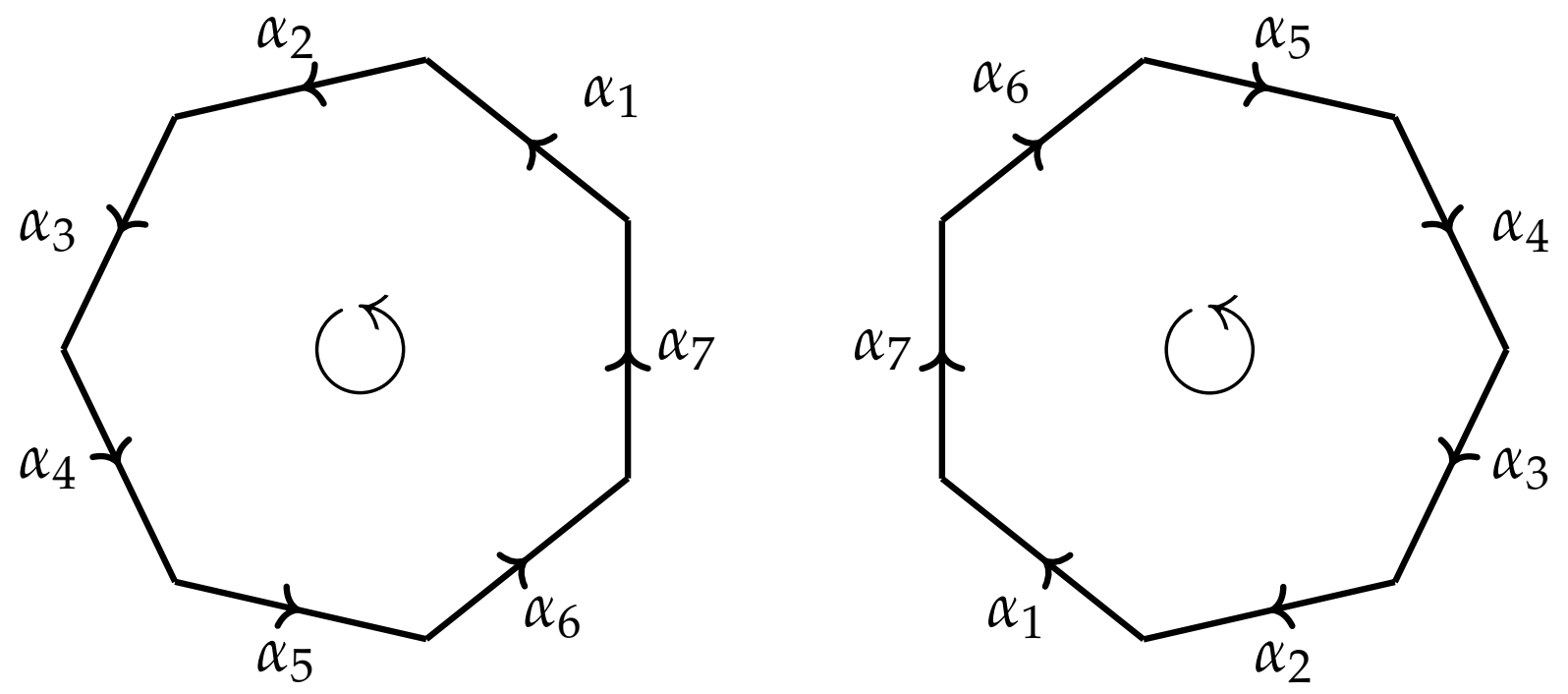

View $\Sigma = \Sigma_g$ as being obtained from two $(2g+1)$-gons, with edges identified as in Figure \ref{fig:4}. Let $\phi: \Sigma \rightarrow \Sigma$ be the homeomorphism obtained by rotating each $2g+1$-gon by $2 \pi /(2g+1)$ counterclockwise about its center. This action is semifree with three fixed points and the cyclic group generated by $\phi$ is $\pi \cong \bb{Z}_{2g+1}$. Let $\alpha_1, \ldots \alpha_{2g+1} \in \pi_1(\Sigma)$ be the generators as illustrated in the figure, and note that these satisfy the relations
$$\alpha_1 \alpha_2 \ldots \alpha_{2g+1} = \alpha_{2g+1} \ldots \alpha_2 \alpha_1 = e.$$

\begin{remark}\label{rem:standard}
(a) The usefulness of this presentation is that the action on the generators is relatively simple. We will see a much more general version of this in Theorem \ref{thm:Action}. 

\medskip

(b) Eliminating $\alpha_{2g+1}$, the remaining variables satisfy the single relation
\begin{equation}\label{eq:standardfunny}
\alpha_1 \alpha_2 \ldots \alpha_{2g} \alpha_{1}^{-1} \alpha_2^{-2} \ldots \alpha_{2g}^{-1} = e.
\end{equation}
These remaining generators $\alpha_1, \ldots, \alpha_{2g}$ can be relabeled to yield a standard generating set with (\ref{eq:standardfunny}) recovering the product of commutators condition. 

\medskip

(c) Let $\psi: \Sigma \rightarrow \Sigma$ be a hyperelliptic involution. The composition $\psi \circ \phi$ is a map on $\Sigma$ of order $4g+2$ which, by Wiman's theorem \cite{Wiman}, is the maximal order a periodic map on $\Sigma = \Sigma_g$ can have. \hfill $\diamond$
\end{remark}

The action of $\phi$ on $\pi_1(\Sigma_g)$ is $\phi_* \alpha_j =\alpha_{j+1} $ for $j \in \bb{Z}_{2g+1}$. Then $\fix_\R(\phi)$ can be identified with the set of tuples $(A_1, \ldots , A_{2g+1})$ in $S^3$ with 
\begin{equation}\label{eq:2g+1}
A_1 A_2 \ldots A_{2g+1} = A_{2g+1} \ldots A_2 A_1 = 1
\end{equation}
and $A_j = A_{j+1}$ for $j \in \bb{Z}_{2g+1}$. Projecting to $A_1$, we thus get a homeomorphism
$$\fix_\R(\phi) \cong \zeta_{2g+1}$$
with the $(2g+1)$-roots of unity in $S^3$.

Now suppose that $N \rightarrow \Sigma$ is a $\bb{Z}_{2g+1}$-equivariant oriented circle bundle with Euler class $k$. In light of Remark \ref{rem:standard}, the fixed points on $\R(N)$ can be identified with those tuples $(A_1, \ldots, A_{2g+1}, Z)$ with 
$$Z^k = A_1 A_2 \ldots A_{2g} A_1^{-1} A_2^{-1} \ldots A_{2g}^{-1},$$
$[A_j, Z] = 1$, and $A_{j+1} = A_j$ for all $j$. Projecting to $(A_1, Z)$, we see that the fixed point set $\fix_\R(\widetilde{\phi}) \subseteq \R(N)$ can be identified with the set of commuting pairs in $\zeta_{2g+1} \times \zeta_{k}$. When either of these pairs is central, then the other components can be any element of $S^3$, so the associated conjugacy class is an interval; otherwise the conjugacy class is a point. Using Proposition \ref{prop:SF}, we obtain the following result about the quotient space $N / \bb{Z}_{2g+1}$.

\begin{cor}
Assume $\bb{Z}_{2g+1}$ acts freely on $N$. The character variety $\X(N^k_g/ \bb{Z}_{2g+1})$ is homeomorphic to a disjoint union of an interval and $g(k-1)$ (resp. $g(k-2)$) points when $k$ is odd (resp. even).  
\end{cor}

We now seek to understand the fixed point set on the character variety $\X(\Sigma)$. We will first determine $\R(\Sigma_\phi)$, and then deduce information about $\fix_\X(\phi)$ and the fixed point set on $\X(N)$ for an equivariant bundle. 

Now $\R(\Sigma_\phi)$ is the set of tuples $(A_1, \ldots , A_{2g+1}, T)$ satisfying (\ref{eq:2g+1}) and $A_{j+1} = T^{-1} A_j T$ for all $j \in \bb{Z}_{2g+1}$. Any such tuple is entirely determined by $(A_1, T)$ and so projecting to these we see that $\R(\Sigma_\phi)$ can be identified with the set $\mathcal{C}$ of pairs $(A, T)$ so that
$$(T A^{-1})^{2g+1} = T^{2g+1}  = (TA)^{2g+1}$$
and with the property that $T^{2g+1} \in C_{A}$ is in the centralizer of $A$. 

If $A^{2g+1} = 1$ and $T$ commutes with $A$, then $(A, T) \in \mathcal{C}$. For $m$ an integer in $[0 ,g]$, write $\mathcal{C}_m$ for the connected component of $\mathcal{C}$ containing $(e^{2 \pi m \hat{i}/(2g+1)}, 1)$. If $(A, T)$ is a commuting pair in $\mathcal{C}$, then $A^{2g+1} = 1$ and $(A, T) \in \mathcal{C}_m$ for $m$ with $\mathrm{Re}(A) = \cos(2 \pi m / (2g+1))$, so the $\mathcal{C}_m$ contain all reducible elements of $\mathcal{C}$. Below we will (i) show that every element of $\mathcal{C}_m$ is reducible, (ii) show the $\mathcal{C}_m$ are pairwise disjoint, and (iii) determine the remaining connected components of $\mathcal{C}$.

We establish (ii) first: Suppose $[0, 1] \ni t \mapsto (A_t, T_t)$ is a path of elements in $\mathcal{C}$ from $(e^{2 \pi \hat{i} m / (2g+1)}, 1)$ to $(e^{2 \pi \hat{i} n / (2g+1)}, 1)$ for some integers $m, n$ with $0 \leq m, n \leq g$. To show $m = n$, it suffices to show that $A_t^{2g+1} = 1$ for all $t$. Working by contradiction, assume there is some $t \in [0, 1]$ with $A_t^{2g+1} \neq 1$, and let $t_*$ be the infimum of all such $t$. It follows that there is some $\epsilon > 0$ so that if $t \in (t_*, t_* + \epsilon)$, then $T_t^{2g+1} \neq \pm 1$. For such $t$, $T_t$ and $A_t$ necessarily commute (since $T_t^{2g+1}$ and $A_t$ commute), and this implies that $A_t^{2g+1} =1$ contrary to the definition of $t_*$.

Next we will determine the components of $\mathcal{C}$ other than the $\mathcal{C}_m$. If $(A, T) \in \mathcal{C}$ and $T^{2g+1} \neq \pm 1$, then $A$ and $T$ commute and $A^{2g + 1} = 1$. This implies $(A, T) \in \mathcal{C}_m$ for some $m$. We thus seek to understand those $(A, T) \in \mathcal{C}$ with $T^{2g + 1} = \pm 1$. Since $(A, T) \in \mathcal{C}$ if and only if $(A, -T) \in \mathcal{C}$, by replacing $T$ with $-T$, we may assume $T^{2g+1} = 1$; note also that the pairs $(A, T)$ and $(A, -T)$ belong to distinct components of $\mathcal{C}$ when $A$ and $T$ do not commute. As such, we are in search of those noncommuting tuples $(A, T)$ with $(TA)^{2g+1} = T^{2g+1} = (TA^{-1})^{2g+1} = 1$. From these identities, we see that the pair $(T^2 A^{-1} T^{-1},T^2)$ satisfy the hypotheses of Lemma \ref{lem:LiorsCountingLemma}. Accounting for the duplicate $T^{2g+1} = -1$ case, we see from that lemma that $\mathcal{C} \backslash \cup_m \mathcal{C}_m$ consists of $\frac{1}{3}g (g-1)(2g-1)$ components. This establishes (iii). Moreover, we also see from this lemma that each irreducible representation is not in a component containing a reducible representation. Thus, each $\mathcal{C}_m$ only contains reducible representations, as claimed in (i), and so 
$$\mathcal{C}_m = \left\{(A, T) \; \Big| \; \mathrm{Re}(A) = \cos(2 \pi m / (2g+1)), \; AT = TA \right\}.$$
For example, $\mathcal{C}_0 = \left\{1 \right\} \times S^3$ and, when $1 \leq m \leq g$, the quotient $\mathcal{C}_m / S^3$ by diagonal conjugation is homeomorphic to $S^1$. We summarize this as follows. 

\begin{prop}
The representation variety $\R(\Sigma_\phi)$ consists of $\frac{1}{3}(2g + 1)(g^2 - 2g + 3)$ connected components. Of these, $g+1$ consist entirely of reducible representations and the union of these components is equivariantly diffeomorphic to $\cup_{0 \leq m \leq g} \mathcal{C}_m$. Each of the remaining $\frac{1}{3}g (g-1)(2g-1)$ components is equivariantly diffeomorphic to $\SO(3)$ (with $S^3$ acting transitively on $\SO(3)$ through the double cover $S^3 \rightarrow \SO(3)$) and consists entirely of irreducible representations. 
\end{prop}

Of course, the character variety $\X(\Sigma_\phi) = \R(\Sigma_\phi)/S^3$ has the same number of components. As for the fixed point set, there are $g+1$ components of $\fix_\X(\phi)$ containing reducible representations (those associated to the components $\mathcal{C}_m$ from above). All other components of $\fix_\X(\phi)$ consist entirely of (equivalence classes of) representations with centralizer $C_\rho = \left\{ \pm 1 \right\}$, and so on these the map (\ref{eq:chimap}) restricts to a 2-1 cover. This gives:

\begin{cor}
The fixed point set $\fix_\X(\phi) \subseteq \X(\Sigma)$ has $\frac{1}{6}(2g^3 - 3g^2 + 7g + 6)$ connected components. One of these is $\mathcal{C}_0 / S^3 \cong [0, 1]$, $g$ of these are homeomorphic to $\mathcal{C}_m / S^3 \cong S^1$, and the remaining components are all points. 
\end{cor}

\subsection{Rotations (\texorpdfstring{$\pi = \mathbb{Z}_g$}{})}\label{sec:Rotations}

\begin{figure}[ht]
\centering
\resizebox{7cm}{6.5cm}{
\begin{tikzpicture}
\pic (upper) at (0,pi) {handle};
\pic[rotate=120] (bl) at (30:-pi) {handle};
\pic[rotate=-120] (br) at (150:-pi) {handle};
\draw (upper-lpoint) to[out=280,in=20] (bl-rpoint);
\draw (upper-rpoint) to[out=260,in=160] (br-lpoint);
\draw (bl-lpoint) to[out=45,in=135] (br-rpoint); 
\node (x0) at (.8,-.3) [ label={\Huge $x_0$}] {};
\node (alpha1) at (1.8,.8) [ label={\Huge $\alpha_1$}] {};
\node (beta1) at (.4,8) [ label={\Huge $\beta_1$}] {};
\node (alpha2) at (-1.8,.1) [ label={\Huge $\alpha_2$}] {};
\node (beta2) at (-7.65,-4.5) [ label={\Huge $\beta_2$}] {};
\node (alpha3) at (.25,-2.8) [ label={\Huge $\alpha_3$}] {};
\node (beta3) at (7.75,-4.5) [ label={\Huge $\beta_3$}] {};
\node (phi) at (7.5,5) [ label={\Huge $\phi$}] {};
\end{tikzpicture}}
\caption{Illustrated here is the homeomorphism $\phi$ for $g = 3$, together with a standard generating set $\alpha_j, \beta_j$ for the fundamental group. \hfill $\diamond$}
\label{fig:3}
\end{figure}
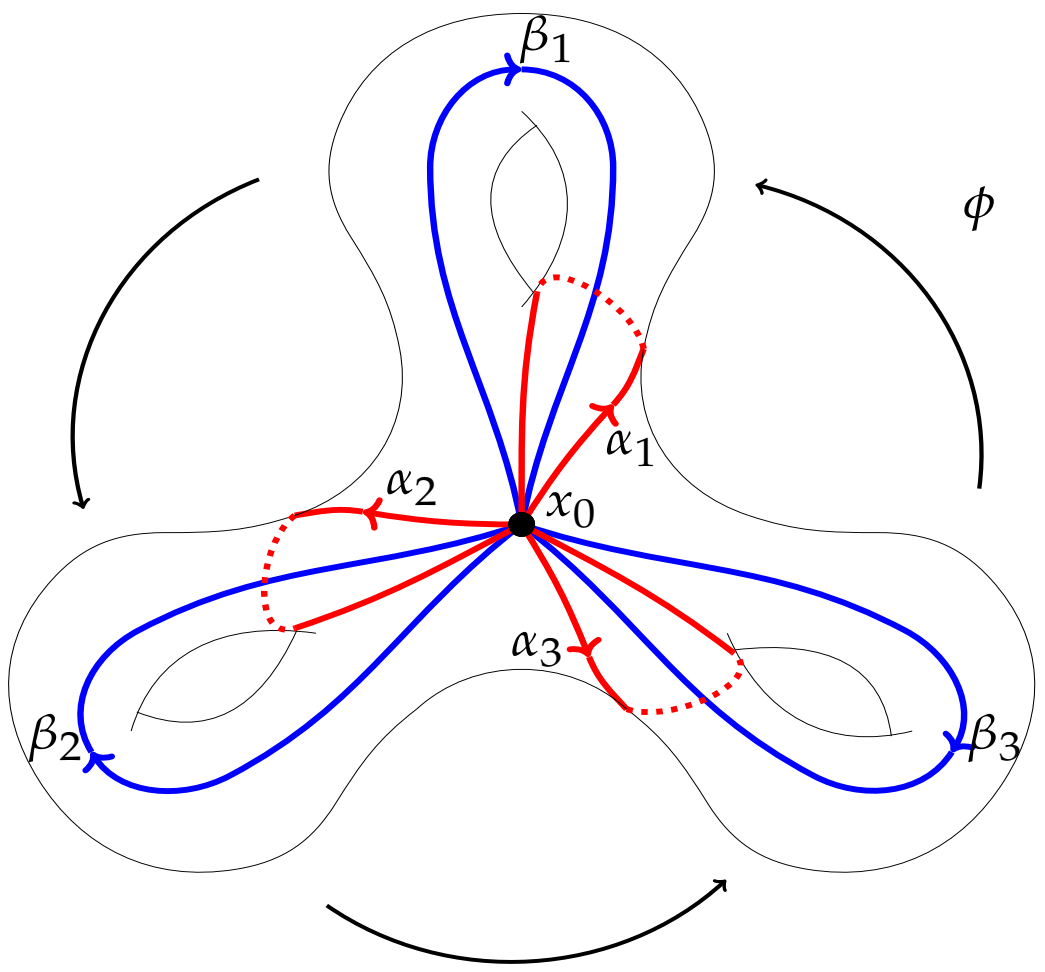

For $g \geq 2$, view $ \Sigma = \Sigma_g$ as being obtained from the 2-sphere by taking the connected sum with $g$ tori, symmetrically spaced along the equator. Rotation by angle $2\pi / g$ about the North--South pole gives a homeomorphism $\phi: \Sigma \rightarrow \Sigma$ with order $g$, exactly two fixed points, and quotient homeomorphic to $T^2$. See Figure \ref{fig:3}.

\begin{prop}\label{prop:Rot1}
The number of connected components of $\X(\Sigma_\phi) $ and of $\R(\Sigma_\phi)$ is $\lfloor g^2/2 \rfloor + 1$.  
\end{prop}

\begin{proof}
The representation variety $\R(\Sigma_\phi)$ can be identified with the set of tuples 
$$((A_j, B_j)_j, T) \in (S^3)^{2g+1}$$ 
satisfying $\prod_{j} [A_j, B_j] = 1$, and $T A_j T^{-1} = A_{j+1}$ and $T B_j T^{-1} = B_{j+1}$ for all $1 \leq j \leq g$. Any such tuple is determined by $T$ together with $A \defeq A_1$ and $B \defeq B_1$. To express the relations entirely in terms of these, first note that $T^g A T^{-g} =A$ and $T^g B T^{-g} = B$, so $T^g \in C_A \cap C_B$ lies in the intersection of the centralizers of $A, B \in S^3$. Also note that $[A_j, B_j] = T^{j-1}[A_1,B_1] T^{-(j-1)}$, so we must have 

\begin{equation} \label{mtrcomm1}
([A, B]T^{-1})^{g-1}  [A,B] T^{g-1} = 1.   
\end{equation}

Conversely, let $\mathcal{C}$ be the set of triples $(A, B, T)$ so that $T^g \in C_A \cap C_B$ and (\ref{mtrcomm1}). Then the above discussion shows there is a homeomorphism $\mathcal{C} \cong \R(\Sigma_\phi)$. Define $\mathcal{C}_0$ to be the component of $\mathcal{C}$ containing the trivial tuple $(1, 1, 1)$.

\medskip

\noindent \emph{Claim 1: The component $\mathcal{C}_0$ contains all tuples with $A, B$ commuting.}

\medskip

If $A, B$ commute, then the condition on $T$ for $(A, B, T) \in \mathcal{C}$ is just that $T^g \in C_A \cap C_B$. Our strategy is to show that we can connect any such tuple $(A, B, T)$ by a path on $\mathcal{C}$ to $(1, 1, 1)$. Since the space of commuting tuples in $S^3$ is connected, it suffices to connect $(A, B, T)$ to $(A, B, 1)$ (by a path in $\mathcal{C}$). We will carry this out under the assumption that $A, B$ are not both central elements of $S^3$, leaving the remaining case to the reader. When $A$ and $B$ are not both central, the intersection $C_A \cap C_B$ is a maximal torus. If $T^g \neq \pm 1$, then $T \in C_A \cap C_B$ and we can find a path $T_t \in C_A \cap C_B$ with $T_0 = T$ and $T_1 = 1$. In this case, $(A, B, T_t)$ is the path we are after. 

If $T^g = \pm 1$, then $(S T S^{-1})^g = \pm 1$ for all $S \in S^3$. Since $C_A \cap C_B$ is a maximal torus, we can find a path $S_t$ with $S_0 = 1$ and $S_1 T S_1^{-1} \in C_A \cap C_B$. Then $(A, B, S_t T S_t^{-1})$ is a path in $\mathcal{C}$ starting at $(A, B, T)$ and ending at a point for which the previous case applies. This finishes the proof of Claim 1. 

\medskip

To determine the remaining connected components, note that when $A, B$ do not commute, the intersection $C_A \cap C_B = \left\{ \pm 1 \right\}$ consists of the center of $S^3$. Then $T^{g} = \pm 1$ and so (\ref{mtrcomm1}) is equivalent to $([A, B] T^{-1})^g = \pm 1$. Consider the function 
$$F: \mathcal{C} \backslash \mathcal{C}_0 \longrightarrow [-1, 1]^2, \indent (A, B, T) \longmapsto \big(\mathrm{Re}(T), \mathrm{Re}([A, B]T^{-1}) \big).$$
The proposition now follows from the next claim.  

\medskip

\noindent \emph{Claim 2: The fibers of $F$ are connected. The image of $F$ consists of the set of pairs $(\cos(\pi k / g), \cos(\pi  \ell / g))$
for distinct integers $k, \ell \in \left\{ 0, 1, \ldots,  g  \right\}$ having the same parity.}

\medskip

Connectivity of the fibers is a consequence of the following three observations: (i) The conjugacy classes in $S^3$ are connected and the conjugacy class of $T$ is determined by $\mathrm{Re}(T)$. (ii) For each $T \in S^3$ and $r \in [-1, 1]$, the set of $z \in S^3$ with $\mathrm{Re}(z T^{-1}) = r$ is connected and nonempty (it is the intersection of $S^3$ with the plane $zT^{-1} + T z^* = 2r$ in $\bb{H}$). (iii) The fibers of the commutator map $(A, B) \mapsto [A, B]$ are connected and nonempty. 

To determine the image of $F$, we have seen that $T$ and $[A, B ] T^{-1}$ are $2g$-roots of unity, so $F(A, B, T) = (\cos(\pi k /g), \cos(\pi \ell / g))$ with $0 \leq k, \ell \leq g$. These integers necessarily have the same parity because $T^{g} = ([A, B] T^{-1})^g = \pm 1$. Suppose $k = \ell$. Then $T^{-1}$ and $[A, B] T^{-1}$ are in the same conjugacy class so, by (ii) and (iii) above, we can path-connect $(A, B, T)$ through $\mathcal{C}$ to an element $(A', B', T)$ with $[A', B'] T^{-1} = T^{-1}$. This implies $[A', B'] = 1$ and so $(A, B, T) \in \mathcal{C}_0$, which is impossible since we have removed $\mathcal{C}_0$ from the domain of $F$. 

Conversely, given distinct $k, \ell$ of the same parity, take $T = e^{\pi \hat{i} k / g}$ and $A, B$ to be any elements with $[A, B] = e^{\pi \hat{i} (k+ \ell)/g}$. Then $(A, B, T) \in \mathcal{C}$ and $F(A, B, T) = (\cos(\pi k/g), \cos(\pi \ell/ g))$. 
\end{proof}

A closer analysis of the proof reveals the homeomorphism type of the nonidentity components of $\R(\Sigma_\phi)$ and $\X(\Sigma_\phi)$. For example, the component containing $(A, B, T)$ with $T \neq \pm 1$ is an $\SO(3)$-bundle over a 2-sphere-bundle over a 2-sphere: The $T$'s with fixed $\mathrm{Re}(T)$ is a 2-sphere; for each such $T$ there is a 2-sphere's worth of $z \in S^3$ with $\mathrm{Re}(z T^{-1})$ fixed (and not equal to $\mathrm{Re}(T)$); for each choice of $z \neq 1$, the set of $(A, B)$ with $[A, B] = z$ is diffeomorphic to $\SO(3)$.

\begin{cor}
The number of connected components of $\fix_\X(\phi) \subseteq \X(\Sigma)$ is $\lfloor (g^2+4)/4 \rfloor$.
\end{cor}

\begin{proof}
As we saw in the proof of Proposition \ref{prop:Rot1}, if $\rho$ is not in the identity component of $\R(\Sigma_\phi)$, then the centralizer of $\rho$ is $\left\{ \pm 1 \right\}$. It follows that (\ref{eq:chimap}) is 2-1 off of the identity component. From Proposition \ref{prop:Rot1} we see that the number of non-identity components of $\fix_\X(\phi)$ is $(g^2-1 )/4$ when $g$ is odd, and $g^2/4$ when $g$ is even. 
\end{proof}

\begin{prop}
The number of connected components of $\fix_\R(\phi) \subseteq \R(\Sigma)$ is $\lfloor g/2 \rfloor + 1$. 
\end{prop}

\begin{proof}
Similar to the proof of Proposition \ref{prop:Rot1}, the fixed point set $\fix_{\R}(\phi)$ can be identified with the set $\mathcal{C}$ of pairs $(A, B)$ with $[A, B]^g = 1$. The components of this fixed point set are detected by the function $F: \mathcal{C} \rightarrow [-1, 1]$ given by $F(A, B) = \mathrm{Re}([A, B])$. The image of $F$ is $\left\{ \cos(2 \pi k / g) \; \vert \; 0 \leq k \leq \lfloor g/2 \rfloor \right\}$, so the proposition follows.  
\end{proof}

\subsection{Symmetrically increasing the genus}\label{sec:Symmetrically}

Fix a $ \pi$-action on $\Sigma_g$ and assume there is a point $x \in \Sigma_g$ with trivial isotropy subgroup, so the orbit $\mathcal{O}_x$ of $x$ has cardinality equal to the order $\mathrm{ord}(\pi)$ of $\pi$. For $h \geq 1$, consider a surface $\Sigma_h$ of genus $h$. Let $S$ be the surface obtained from $\Sigma_g$ by taking the connected sum with $\Sigma_h$ \emph{at each} element of $\mathcal{O}_x$. Then $S$ is a closed, connected, oriented surface of genus $g + h \cdot \mathrm{ord}(\pi)$. The action of $\pi$ on $\Sigma_g$ gives an action of $\pi$ on $S$ with the same fixed point stratum and with quotient $S / \pi$ given by the connected sum of $\Sigma_g / \pi$ with $\Sigma_h$. We will write $\Sigma_g \#_\pi \Sigma_h$ in place of $S$. 

\begin{ex}
Consider the standard action of $\pi = \bb{Z}_g$ on $S^2$ by rotation about the polar axis. Then $S^2 \#_{\pi} \Sigma_1$ is the surface of genus $g$ and the induced $\pi$-action is precisely that of Section \ref{sec:Rotations}. \hfill $\diamond$
\end{ex}

For a more involved example, consider the action of $\pi = \bb{Z}_{2g+1}$ on $\Sigma_g$ from Section \ref{sec:SomeMaps3}. The presentation of $\pi_1(\Sigma_g)$ from that section combines with the standard presentation on each copy of $\pi_1(\Sigma_h)$, to give a presentation of the fundamental group of $S = \Sigma_{g} \#_\pi \Sigma_h$ with generators $\alpha_1, \ldots, \alpha_{2g+1}$ and $\gamma_{1, \ell}, \delta_{1, \ell}, \ldots, $ $\gamma_{ h, \ell}, \delta_{ h, \ell}$ for $1 \leq \ell \leq 2g+1$ subject to the relations
$$\alpha_1 \alpha_2 \ldots \alpha_{2g+1} = \alpha_{2g+1} \ldots \alpha_2 \alpha_1 = 1$$
and 
$$(\alpha_1 \alpha_2 \ldots \alpha_{2g})(\alpha_{2g} \ldots \alpha_2 \alpha_1) \prod_{\ell = 1}^{2g+1} \prod_{i = 1}^h [ \gamma_{i, \ell}, \delta_{i, \ell}] = 1.$$
These generators can be arranged so that $1 \in \bb{Z}_{2g+1}$ acts by
$$\alpha_\ell \longmapsto \alpha_{\ell+1}, \hspace{.5cm} \gamma_{ i, \ell} \longmapsto \gamma_{i, \ell+1}, \hspace{.5cm} \delta_{i, \ell} \longmapsto \delta_{i, \ell+1}$$
for $\ell \in \bb{Z}_{2g+1}$ and $1 \leq i \leq h$. 

Let $\rho \in \R(S)$ be a fixed point under this action. This is uniquely determined by $A \defeq \rho(\alpha_1)$, $C_i \defeq \rho(\gamma_{i, 1})$, and $D_i \defeq \rho(\delta_{i, 1})$ for $1 \leq i \leq h$. The condition on these is just that
$$A^{2g+1} = 1, \hspace{.5cm} \Big( \prod_{i = 1}^{h} [C_i, D_i] \Big)^{2g+1} = 1.$$
From this, we see that the fixed points set $\fix(\bb{Z}_{2g+1}, \R(S))$ is homeomorphic to the Cartesian product of $\zeta_{2g+1}$ and the solutions of 
$$\Big( \prod_{i = 1}^{h} [C_i, D_i] \Big)^{2g+1} = 1.$$
The solutions of this latter equation are described by the \emph{twisted} representation varieties
$$\R_\eta(\Sigma_h) \defeq \left\{(C_i, D_i)_i \; \Bigg| \; \prod_i [C_i, D_i] = \eta \right\}$$
for $\eta \in \zeta_{2g+1}$, the topological type of which is fairly well-understood (e.g., it is connected). As a concrete consequence, we have the following.

\begin{cor}
The fixed point set $\fix(\pi, \R(S))$ has $(g+1)^2$ connected components. 
\end{cor}

\section{Semifree cyclic actions, in general}\label{sec:InGeneral}

Throughout this section we fix an integer $n \geq 2$, and assume $\Sigma$ admits an effective semifree action of $\bb{Z}_n$ through orientation-preserving homeomorphisms. Then the fixed points are necessarily isolated; enumerate them as $x_1, \ldots, x_t$. Write $k_j \in \bb{Z}_n$ for the rotation number of $x_j$, and we recall that this is necessarily a unit in $\bb{Z}_n$; write $k^{*}_j \in \bb{Z}_n$ for the multiplicative inverse of $k_j$. It is necessarily the case that the sum $\sum_{j} k^{*}_j \equiv 0 \in \bb{Z}_n$ vanishes mod $n$ \cite{Nielson}, \cite[Prop. 1.1]{EE}; we refer to $\sum_{j} k^{*}_j \equiv 0$ as the \emph{consistency condition}. To simplify the discussion, we have assumed that the quotient $\Sigma / \bb{Z}_n$ is a topological 2-sphere, which is equivalent to assuming $g = (t - 2)(n-1)/2$; see Remark \ref{rem:genus} for an approach to handling higher genus quotients. 

In this section we give a description of the fixed $\bb{Z}_n$-fixed point set on $\R(\Sigma)$. As a starting point we take the following result giving a favorable presentation for $\pi_1(\Sigma)$ that is symmetric relative to the action of $\bb{Z}_n$. To state it, for $2 \leq j \leq t$ and $1 \leq i \leq k_j^*$ we define
\begin{equation}\label{eq:elldef}
\begin{array}{rcl}
\ell(i, j) & \defeq & k_1( i - k^{*}_j  - k^{*}_{j+1} - k^{*}_{j+2} - \ldots - k^{*}_t) \\
&\equiv & k_1(i + k_1^* + \ldots + k_{j-1}^*) \pmod n.
\end{array}
\end{equation}

\begin{thm}\label{thm:Action}
Suppose $\Sigma$ is a surface with genus $g= (t - 2)(n-1)/2$ admitting a semifree and orientation-preserving action of $\bb{Z}_n$ with weights $k_1, \ldots, k_t$. Let $x_1$ be the fixed point of this action having weight $k_1$. Then the fundamental group $\pi_1(\Sigma, x_1)$ has a presentation with generators $\alpha_{j, \ell}$ for $2 \leq j \leq t$ and $\ell \in \bb{Z}_n$, subject to the relations $\prod_{i = 0}^{n-1} \alpha_{j, i k_1} = e$ for all $2 \leq j \leq t $ and
$$\prod_{j = 2}^t \prod_{i = 1}^{k^{*}_j} \alpha_{j, \ell(i, j)+ \mu} = e$$
for $\mu \in \bb{Z}_n$, with $\ell(i, j)$ as in (\ref{eq:elldef}). These generators can be chosen so that the action of $s \in \bb{Z}_n$ on $\pi_1(\Sigma, x_1)$ sends $\alpha_{j, \ell} \mapsto \alpha_{j, \ell - k_1 s}$. 
\end{thm}

We have relegated the proof of Theorem \ref{thm:Action} to Section \ref{sec:ProofOfTheoremAction}. 

\begin{remark}\label{rem:genus}
(a) Our proof will be constructive and will show that if $k_1, \ldots, k_t$ are any invertible elements of $\bb{Z}_n$ satisfying the consistency condition, then there is a semifree action of $\bb{Z}_n$ on $\Sigma_{g}$ with $t$ fixed points, $g = (t - 2)(n-1)/2$, and weights given by the $k_j$. By Nielson's theorem \cite{Nielson}, this action is unique up to conjugation in $\homeo(\Sigma_g)$. 

The genus $g$ can be arranged to be larger than $(t - 2)(n-1)/2$ using the strategy of Section \ref{sec:Symmetrically}. By Nielson's theorem \cite{Nielson}, we can, in this way, obtain \emph{every} semifree and orientation-preserving action of $\bb{Z}_n$ on any surface $\Sigma$ (closed, connected, and oriented, as usual). 

\medskip

(b) The relations of Theorem \ref{thm:Action} can be conceptualized as follows: The quotient $\Sigma / \bb{Z}_n$ is a topological 2-sphere with $t$ singular (orbifold) points $[x_1], \ldots , [x_t]$. Use $[x_1]$ as a basepoint. For $2 \leq j \leq t$, let $\gamma_j$ be a small loop in this quotient starting at $[x_1]$ and circling $[x_j]$; we can arrange the $\gamma_j$ so they intersect only at the basepoint and so the ``inside'' of $\gamma_j$ (the part containing $[x_j]$) is a topological disk. The preimage in $\Sigma_g$ of the inside of $\gamma_j$ is a topological disk in $\Sigma_g$ containing $x_j$; we can view this disk as giving rise to the relation $\prod_i \alpha_{j, i k_1} = e$. Since $\Sigma_g / \bb{Z}_n$ is a 2-sphere, the complement of the insides of the $\gamma_j$ is also a topological disk. This disk contains no fixed points, so its preimage in $\Sigma_g$ is a disjoint union of $n$ topological disks, and we can view these as giving rise to the relations $\prod_j \prod_i \alpha_{j, \ell(i, j) + \mu} = e$ for $\mu \in \bb{Z}_n$. \hfill $\diamond$
\end{remark}

\begin{ex}
Consider the case where $n = 2$. Then $t$ is necessarily even and the surface $\Sigma$ of Theorem \ref{thm:Action} has genus $g= t/2 - 1$. Since $n = 2$, the only possible weight is $k_j = 1$ for all $j$. Since $t$ is even the consistency condition $\sum_j k^{*}_j \equiv 0$ is automatic. By Nielson's classification, this is just the hyperelliptic involution, which we considered in Section \ref{sec:HEI}. 

To better-understand the presentation given by Theorem \ref{thm:Action}, consider the case where $t = 4$ (so $g = 1$). The presentation  has six generators $\alpha_{j, \ell}$ with $2 \leq j \leq 4$ and $\ell \in \bb{Z}_2$. The index $i$ in $\ell(i, j)$ can only be $1$, and a computation shows $\ell(i, j) = \ell(1, j) \equiv j \pmod{2}$ for all $j$. Using the relations of Theorem \ref{thm:Action}, we see that $\pi_1(\Sigma)$ is freely generated, as an abelian group, by $\alpha_{2, 0}$ and $\alpha_{3, 1}$, as expected given that we know $\Sigma$ must be a torus.  \hfill $\diamond$
\end{ex}

\begin{ex}
 When $t = n = 3$, Theorem \ref{thm:Action} gives a presentation with 6 generators $\alpha_{j, \ell}$ for $j = 2, 3$ and $\ell \in \bb{Z}_3$. Then $\ell(i, j) = \ell(1, j) \equiv j \pmod{3}$, from which we see that $\pi_1(\Sigma)$ is freely generated, as an abelian group, by $\alpha_{2, 0}$ and $\alpha_{2, 1}$. Of course, this is expected since the genus is $g = (t - 2)(n-1)/2 = 1$. \hfill $\diamond$
\end{ex}

With these fundamental group preliminaries, we now shift to the representation varieties. As usual, each $\rho \in \R(\Sigma)$ is determined by the values $A_{j, \ell} \defeq \rho(\alpha_{j, \ell})$ for $2 \leq j \leq t$ and $\ell \in \bb{Z}_n$. The fixed point condition for $\fix_\R(\bb{Z}_n)$ implies that
$$A_{j, 0} = A_{j, 1} = \ldots = A_{j, n-1}$$
for each $2 \leq j \leq t$, so we only need to keep track of $A_j \defeq A_{j, 0}$. Then we see that $\rho \in \R(\Sigma)$ is a fixed point if and only if
$$A_2^n = \ldots = A_t^n = A_2^{k_2^*} A_3^{k_3^*} \ldots A_t^{k_t^*} = 1.$$
Since each of the $k_j^*$ is relatively prime to $n$, we have $A_j \in \zeta_n$ if and only if $B_{j+1} \defeq A_j^{k_j^*} \in \zeta_n$. Then we are looking for those tuples $B_1, \ldots, B_{t-1} \in \zeta_n$ with $B_1 B_2 \ldots B_{t-1} = 1$. Using this to solve for $B_{t-1}$, we obtain the following.

\begin{prop}\label{prop:RepGen}
The fixed point set $\fix_\R(\bb{Z}_n) \subseteq \R(\Sigma)$ is homeomorphic to the set of $(B_1, B_2, \ldots, B_{t-2}) \in \zeta_n^{t-2}$ so that the product $B_1 B_2 \ldots B_{t-2} \in \zeta_n$ is an $n$th root of unity. 
\end{prop}

We thus arrive at the unexpected (for us, anyway) observation that the homeomorphism type of the fixed point set is independent of the weights $k_1, \ldots, k_t$.

Now suppose $N = N^k \rightarrow \Sigma$ is a $\bb{Z}_n$-equivariant oriented circle bundle. Via the clutching construction about a small neighborhood of the fixed point $x_t \in \Sigma$, we find that the fundamental group of $N$ has a presentation with generators $z$ and $\alpha_{j, \ell}$ with $2 \leq j \leq t$ and $\ell \in \bb{Z}_n$ with the property that $z$ is central and the $\alpha_{j, \ell}$ satisfy the same relations from Theorem \ref{thm:Action}, except that the condition $\alpha_{t, 0} \alpha_{t, 1} \ldots \alpha_{t, n-1} = e$ should be replaced by $\alpha_{t, 0} \alpha_{t, 1} \ldots \alpha_{t, n-1} = z^{-k}$. At this point we find it convenient to assume that $k_t = 1$, which can always be arranged by selecting an appropriate generator of $\bb{Z}_n$ (recall the definition of the rotation numbers depends on a choice of generator). Then a straight-forward extension of the above argument shows that the fixed point set on $\R(N)$ can be identified with the tuples $A_2, \ldots, A_{t}, Z$ with $Z$ commuting with all $A_j$, $A_2^n = \ldots = A_{t-1}^n = 1$, $A_t^n = Z^{-k}$, and 
$$A_1^{k_1^*}  \ldots A_{t-1}^{k_{t-1}^*} A_{t} = 1.$$ 
Again setting $B_{j} = A_{j-1}^{k_j^*}$, we have the following.

\begin{prop}\label{prop:RepNGen}
The fixed point set $\fix(\bb{Z}_n, \R(N)) \subseteq \R(N)$ is homeomorphic to the set of tuples $(B_1, \ldots, B_{t-2}, Z) \in \zeta_n^{t-2} \times S^3$ so that $[B_j , Z] = 1$ for all $1 \leq j \leq t-2$ and
$$(B_1 B_2 \ldots B_{t-2})^n = Z^k.$$
\end{prop}

In Section \ref{sec:UpTo4} we compute these fixed point sets for $t\leq 4$. In Section \ref{sec:t=5} we extend this further to $t = 5$ and $n = 3$, and we indicate the limitations of our approach for larger values of $t$.

\subsection{Up to four fixed points} \label{sec:UpTo4}

We continue to assume that there is at least one fixed point, and there can be no action with exactly one fixed point. When $t = 2$ the space $\Sigma = S^2$ is the 2-sphere, the representation variety for which consists of a point. 

As such, the first case that is interesting and relevant is $t = 3$. By Proposition \ref{prop:RepGen}, the fixed point set $\fix_\R(\bb{Z}_n) \subseteq \R(\Sigma)$ is homeomorphic to the space $\zeta_n$ of $n$th roots of unity, and by Proposition \ref{prop:RepNGen}, the fixed point set $\fix(\bb{Z}_n, \R(N))$ is homeomorphic to the set of commuting pairs in $\zeta_n \times \zeta_k$. Alternatively, we can use the general observation that these fixed point spaces are independent of the weights, and just appeal to the analysis of Section \ref{sec:SomeMaps3}, which corresponds to the case $k_1 = k_2 = 1$ and $k_3 = n-2$.

When $t = 4$, Proposition \ref{prop:RepGen} gives that $\fix_\R(\bb{Z}_n)$ is homeomorphic to the space pairs $(A, B) \in \zeta_n \times \zeta_n$ with $AB \in \zeta_n$. For odd $n$, the topology of this space is just as described in Lemma \ref{lem:LiorsCountingLemma}. We leave to the reader the analogous description for even $n$ (here the central element $-1$ plays a distinguished role).

Now suppose $N \rightarrow \Sigma$ has Euler number $k \neq 0$. The following result can be obtained by an argument similar to the one given for Lemma \ref{lem:LiorsCountingLemma}.

\begin{cor}
Suppose $n$ is odd. For even $k$, the fixed point set $\fix(\bb{Z}_n, \R(N))$ is $S^3$-equivariantly homeomorphic to
\[\{*\} \sqcup \{*\} \sqcup \left(\bigsqcup_{i = 1}^{(k n^2 - 2)/2} S^2\right) \sqcup \left(\bigsqcup_{i = 1}^{(n^3 - 6n^2 + 11 n - 6)/12} \SO(3)\right),\]
where $*$ is a point. For odd $k$, the fixed point set $\fix(\bb{Z}_n, \R(N))$ is $S^3$-equivariantly homeomorphic to
\[\{*\} \sqcup \left(\bigsqcup_{i = 1}^{(k n^2-1)/2} S^2\right) \sqcup \left(\bigsqcup_{i = 1}^{(n^3 + 3 n^2 + 5n - 9)/12} \SO(3) \right).\]
\end{cor}

\subsection{Five fixed points with \texorpdfstring{$\bb{Z}_3$-action}{} }\label{sec:t=5}

When $t = 5$, we see from Proposition \ref{prop:RepGen} that $\fix_\R(\bb{Z}_3) \subseteq \R(\Sigma)$ is homeomorphic to the set $\mathcal{C}$ of tuples $(A, B, C, D) \in \zeta_3^4$ satisfying $AB = CD$. To describe this, consider the map $\Pi:\zeta_3\times \zeta_3 \to \bb{R}^3$ defined by $\Pi(A, B) = AB$. Our space $\mathcal{C}$ is the fiber product of $\Pi$ with itself. We will first analyze the fibers of $\Pi$, and then describe the fiber product.  

The image of $\Pi$ is the closed ball $D \subseteq S^3$ of radius $2 \pi / 3$ centered at $1 \in S^3$, the boundary of which is the non-identity component of $\zeta_3$:
$$\zeta_3 = \left\{1 \right\} \sqcup \partial D.$$
The fiber $\Pi^{-1}(A)$ of any $A \in \partial D$ is $\left\{(1, A), (A, 1), (A^{-1}, A^{-1})\right\}$. As such, $\Pi^{-1}(\partial D)$ is the union of three 2-spheres
$$S_1 \defeq \left\{(1, A) \; \big| A \in \partial D \right\}, \hspace{.5cm} S_2 \defeq \left\{(A, 1) \; \big| A \in \partial D\right\},$$
$$S_3 \defeq \left\{(A, A) \; \big| A \in \partial D \right\}.$$
The fiber $\Pi^{-1}(0)$ is a union of the point $\left\{(1, 1) \right\}$ and the 2-sphere
$$S_4 = \left\{(A, A^{-1}) \; \big| A \in \partial D \right\}.$$
Of course, $S_3$ and $S_4$ intersect. The preimage of the intermediary cylinder $\mathrm{int}(D) \backslash \left\{1 \right\}$ is connected and contains $S_3$ and $S_4$ in its closure. Write $W$ for the connected component of $\zeta_3 \times \zeta_3$ containing this latter preimage (and thus also containing $S_3$ and $S_4$). Taking the fiber product, in $\mathcal{C}$ we have a copy of $W \times_\Pi W$, which is the fiber product of the restriction $\Pi \vert_{W}$ with itself. This is connected and contains $S_3 \times S_3$ as well as $S_4 \times S_4$. The other components are $\left\{(1, 1, 1, 1)\right\}$, $\left\{(1, 1) \right\} \times S_4$, $S_4 \times \left\{(1, 1) \right\}$ and $S_i \times S_j$ for $1 \leq i, j \leq 3$ not both equal to 3. From this we conclude that $\fix_\R(\bb{Z}_3) \subseteq \R(\Sigma)$, is equivariantly homeomorphic to 
\[\{*\} \sqcup S^2 \sqcup S^2 \sqcup \left(\bigsqcup_{i = 1}^8(S^2 \times S^2)\right) \sqcup \Big(W \times_\Pi W\Big).\]

\section{A symmetric presentation --- Proof of Theorem \ref{thm:Action}}\label{sec:ProofOfTheoremAction}

Given weights $k_1, \ldots, k_t$ satisfying the consistency condition, we will construct a surface $\Sigma$ equipped with a $\bb{Z}_n$-action having precisely these weights and with $\pi_1(\Sigma)$ having the claimed presentation. By Nielson's classification \cite{Nielson}, any other $\bb{Z}_n$-action is conjugate to this one by an orientation-preserving homeomorphism of $\Sigma$, so Theorem \ref{thm:Action} will follow.

\begin{figure}[ht]
    \centering
\begin{tikzpicture}[
dot/.style={draw,fill,circle,inner sep=1pt}
]
\draw (0,0) ellipse (3*.8cm and 1.5*.8cm) ;
\draw plot [smooth, tension=1.5] coordinates { (1.65*.8,-.075*.8) (0,.6*.8) (-1.65*.8,-.075*.8)};
\draw plot [smooth, tension=1.5] coordinates { (1.8*.8,0.15*.8) (0,-.6*.8) (-1.8*.8,0.15*.8)};
\draw (-2.2,0) ellipse (.15cm and .12cm) ;
\draw (-1.6,0) ellipse (.15cm and .12cm) ;
\draw (1.6,0) ellipse (.15cm and .12cm) ;
\draw (2.2,0) ellipse (.15cm and .12cm) ;
 \node[above left] at (-2, -1.3) {$\Sigma_\circ$};
 \node[left] at (0, -2) {$\Pi_\circ$};
\draw [->] (2.4,.8) arc (120:60:1.2cm);
\draw [-to](0,-1.5) -- (0,-2.5);
\draw (6,0) ellipse (3*.8cm and 1.5*.8cm) ;
\draw plot [smooth, tension=1.5] coordinates { (6+1.65*.8,.075*.8) (6+0,.6*.8) (6-1.65*.8,.075*.8)};
\draw plot [smooth, tension=1.5] coordinates { (6+1.8*.8,0.15*.8) (6+0,-.6*.8) (6-1.8*.8,0.15*.8)};
\fill (6-2.2,0) circle (.05);
\fill (6-1.6,0) circle (.05);
\fill (6+2.2,0) circle (.05);
\fill (6+1.6,0) circle (.05);
 \node[above right] at (6+2, -1.3) {$\Sigma$};
 \draw [-to](6,-1.5) -- (6,-2.5);
\draw plot [smooth, tension=1.5] coordinates { (2*1.65*.8,-2*.075*.8-4) (0,-2.7) (-2*1.65*.8,-2*.075*.8-4)};
\draw plot [smooth, tension=1.5] coordinates { (1.65*.8,-.075*.8-4) (0,.6*.8-4) (-1.65*.8,-.075*.8-4)};
\draw (-2.485,-4.14) ellipse (.15cm and .12cm) ;
\draw (-1.47,-4.14) ellipse (.15cm and .12cm) ;
\draw (2.485,-4.14) ellipse (.15cm and .12cm) ;
\draw (1.47,-4.14) ellipse (.15cm and .12cm) ;
\draw plot [smooth, tension=1.5] coordinates { (-2.335,-4.14) (-1.9775, -4)  (-1.62,-4.14)};
\draw plot [smooth, tension=1.5] coordinates { (2.335,-4.14) (1.9775, -4)  (1.62,-4.14)};
 \node[above left] at (-2, -3.3) {$S$};
  \node[below right] at(-2.485,-4.14){$\gamma_1$};
    \node[below right] at (-1.47,-4.14){$\gamma_2$};
      \node[below left] at(2.485,-4.14){$\gamma_4$};
    \node[below left] at (1.47,-4.14){$\gamma_3$};
\draw [->] (2.4,-3.2) arc (120:60:1.2cm);
\draw plot [smooth, tension=1.5] coordinates { (6+2*1.65*.8,-2*.075*.8-4) (6,-2.7) (6-2*1.65*.8,-2*.075*.8-4)}; 
\draw plot [smooth, tension=1.5] coordinates { (6+1.65*.8,-.075*.8-4) (6,.6*.8-4) (6-1.65*.8,-.075*.8-4)}; 
\draw plot [smooth, tension=1.5] coordinates { (6-1.65*.8,-.075*.8-4) (6-1.9775, -3.9) (6-2*1.65*.8,-2*.075*.8-4) }; 
\draw plot [smooth, tension=1.5] coordinates { (6+1.65*.8,-.075*.8-4) (6+1.9775, -3.9)  (6+2*1.65*.8,-2*.075*.8-4)  }; 
 \node[above right] at (6+2, -3.3) {$S^2$};
 \end{tikzpicture}
  \caption{Illustrated here is the case $n = 2$ and $t = 4$. The vertical arrows are $\mathbb{Z}_n$-quotient maps, and the horizontal arrows are inclusions. Each loop $\gamma_j$ in $S$ encircles one of the $t$ holes appearing in $S$. The solid dots in $\Sigma$ represent the fixed points. The bottom right surface is homeomorphic to $S^2$ with $t$ orbifold points associated to the fixed points in $\Sigma$. \hfill $\diamond$}
\label{fig:10}
\end{figure}
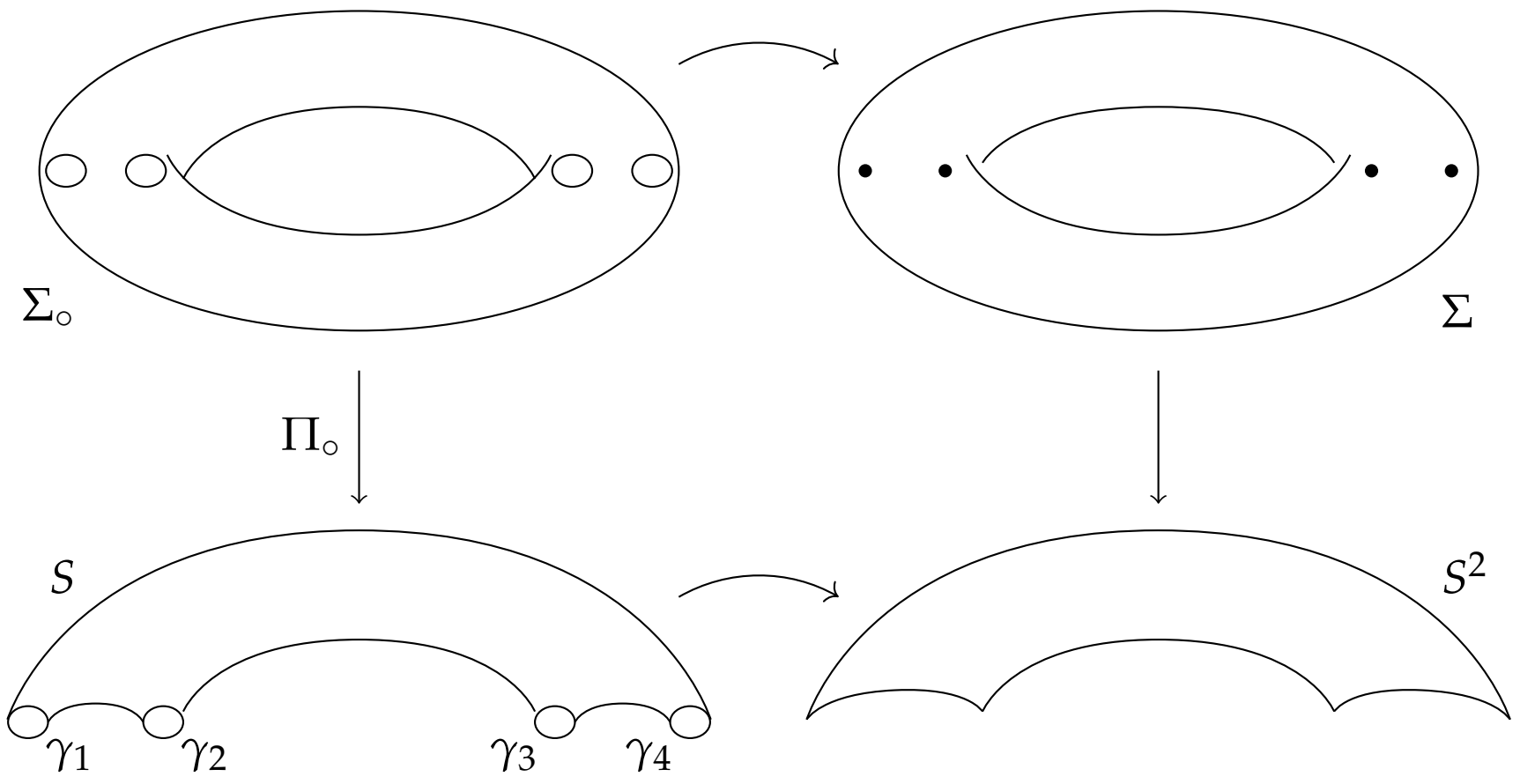

To construct $\Sigma$, remove from $S^2$ small disjoint disks $D_1, \ldots, D_t \subseteq S^2$ to get a surface 
$$S = S^2 - \sqcup^t_{j = 1} \mathrm{int}(D_j)$$ 
with $t$ boundary components. Fix a basepoint $x \in S$ and assume that the disks are arranged counterclockwise around $x$ in the order $D_1$, then $D_2$, then $D_3$, etc. Let $\gamma_j$ be a path that emanates from this basepoint to the boundary of $D_j$, traverses this boundary counterclockwise, and then returns back to $x$ the way it came; this can be arranged so the $\gamma_j$ do not pairwise intersect except at the basepoint. See the bottom left in Figure \ref{fig:10}.

The product $\gamma_1 \gamma_2 \ldots \gamma_t$ is contractible relative to the basepoint $x$. We may therefore view $\pi_1(S, x)$ as being freely generated by $\gamma_j$ for $1 \leq j \leq t -1$, and we set $\gamma_t \defeq (\gamma_1 \gamma_2 \ldots \gamma_{t-1})^{-1}$. Let $N$ be the normal subgroup of $\pi_1(S, x)$ that is generated by the words $\gamma_j^n$ and $\gamma_j^{-k_j} \gamma_1^{k-1}$ for $1 \leq j \leq t - 1$. The next lemma follows from a direct computation and the consistency condition.

\begin{lemma}\label{lem:NN'}
The normal subgroup $N$ contains $\gamma_t^n$ and $\gamma_t^{-k_t} \gamma_1^{k_1}$.
\end{lemma}

The next proposition articulates the key property that $N$ satisfies. 

\begin{prop}\label{prop:quot}
The quotient $\pi_1(S, x) / N$ is a cyclic group of order $n$. For any $1 \leq j \leq t$ and any $a \in \bb{Z}$ coprime to $n$, the coset $\gamma_j^a N$ is a generator of $\pi_1(S, x) / N$. 
\end{prop}

\begin{proof}
We will show that the coset $\gamma_1 N$ has order $n$ and generates the quotient; a similar argument shows that $\gamma_j N$ does the same thing and, since $a$ is coprime with $n$, it follows that each $\gamma_j^a N$ generates the quotient as well.

Let $K \trianglelefteq N$ be the normal subgroup of $N$ generated by all of the generators of $N$ except for $\gamma_1^n$. Then $K$ is a proper subgroup of $N$ since it does not contain $\gamma_1^n$. Thus, the quotient $N / K$ is not trivial. Indeed, $N / K \cong \langle \gamma_1^n \rangle$ is isomorphic to the group freely generated by $\gamma_1^n$. 

We claim that $\gamma_1 K$ generates $\pi_1(S, x) / K$. Indeed, since $\gamma_j^{-k_j} \gamma_1^{k_1} \in K$, it follows that
$$\gamma_j^{k_j} K  = \gamma_1^{k_1} K.$$
We have 
$$\gamma_j K = (\gamma_j^{k_j})^{k^{*}_j} K = (\gamma_1^{k_1})^{k^{*}_j} K = (\gamma_1 K)^{k_1 k^{*}_j}.$$
Thus, the coset of each generator in $\pi_1(S, x)$ can be expressed as a power of $\gamma_1 K$ and so this proves the claim. We also see from this that when all other $\gamma_j $ are expressed in terms of $\gamma_1$, there are no other relations in $K$, so 
$$\pi_1(S, x) / K  = \langle \gamma_1 K \rangle$$
is freely generated by $\gamma_1 K$. Since $N/ K = \langle \gamma_1^n K \rangle$ is freely generated by $\gamma_1^n K$, it follows that
$$\frac{\pi_1(S, x) }{ N} \cong \frac{\pi_1(S, x) / K}{N/ K} \cong \frac{\langle \gamma_1 K \rangle}{\langle \gamma_1^n K \rangle} \cong \frac{\bb{Z}}{n \bb{Z}}.$$
\end{proof}

We will view $\gamma_1^{k_1} N$ as the \emph{preferred} generator of $\pi_1(S, x) / N$; as we will see in Lemma \ref{lem:preferredgen}, this choice will ensure that our construction returns the specified rotation numbers $k_j$. In particular, this choice of generator allows us to identify $\pi_1(S, x) / N$ with $\bb{Z}_n$ via the isomorphism

\begin{equation}\label{eq:PrefIso}
\bb{Z}_n \longrightarrow \pi_1(S, x)/N, \hspace{.5cm} a \longmapsto \gamma_1^{a k_1}.
\end{equation}

Let $\widetilde{S}$ be the universal cover of $S$, and recall that $\pi_1(S, x)$ acts freely on $\widetilde{S}$ with quotient $S$. Define 
$$\Sigma_\circ \defeq \widetilde{S} / N$$ 
to be the quotient of $\widetilde{S}$ by the subaction of $N \trianglelefteq  \pi_1(S, x)$. Since this action is free, this quotient $\Sigma_\circ$ is naturally a connected, oriented surface with boundary. The group $\bb{Z}_n \cong \pi_1(S , x) / N$ acts freely on $\Sigma_\circ$ and the quotient
$$\frac{\Sigma_\circ }{ \bb{Z}_n} \cong  \frac{\widetilde{S}  / N}{\pi_1(S , x) / N} \cong \frac{\widetilde{S} }{\pi_1(S, x)} \cong S$$
recovers our holey sphere. Through this, the projection $\Sigma_\circ \rightarrow \Sigma_\circ / \bb{Z}_n$ gives a surjection
$$\Pi_\circ : \Sigma_\circ \longrightarrow S.$$
Since $\widetilde{S}$ is simply-connected, we have $\pi_1(\Sigma_\circ, x_\circ) \cong N,$ where $x_\circ$ is any element with $\Pi_\circ(x_\circ ) = x$. Since $\bb{Z}_n$ is a finite group, the surface $\Sigma_\circ$ is compact.

\begin{lemma}\label{lem:preferredgen}
The surface $\Sigma_\circ$ has $t$ boundary components. Moreover, the projection $\Pi_\circ$ restricts to each boundary circle of $\Sigma_\circ$ to be an $n$-to-1 map onto a boundary component of $S$, and this induces a bijection $\pi_0(\partial \Sigma_\circ) \cong \pi_0(\partial S)$. In terms of the action of $\bb{Z}_n \cong \pi_1(S, x) / N$ on $\Sigma$, the action of the generator $\gamma_1^{k_1} N$ on the $j$th boundary component of $\Sigma_\circ$ is by rotation by $e^{2 \pi \hat{i} k_j/n}$.
\end{lemma}

\begin{proof}
View $\widetilde{S}$ as the based-homotopy classes $[\gamma]$ of paths $\gamma: [0, 1 ]\rightarrow  S$ with $\gamma(0) = x$, so the covering map $\widetilde{S} \rightarrow S$ is given by $[\gamma] \mapsto \gamma(1)$. The elements of the boundary of $\widetilde{S}$ are precisely the homotopy classes of paths that terminate in $\partial S$. Any one component of $\partial \widetilde{S}$ lying over $\partial D_j$ can be viewed as the universal cover of $\partial D_j \cong S^1$, and the infinite cyclic subgroup $\langle \gamma_j \rangle \leq \pi_1(S, x)$ is precisely the group of deck transformations. Moreover, this cyclic subgroup $\langle \gamma_j \rangle$ is the subgroup that preserves the components lying over $\partial D_j$: If $\widetilde{y} \in \partial \widetilde{S}$ descends to $\partial D_j \subseteq \partial S$ under the projection, and $\gamma \in \pi_1(S, x)$, then $\gamma \cdot \widetilde{y}$ and $\widetilde{y}$ lie in the {same} component of $\partial \widetilde{S}$ if and only if $\gamma = \gamma_j^k$ for some $k \in \bb{Z}$. We also have that $\gamma \cdot \widetilde{y}$ and $\widetilde{y}$ descend to the same element of $\widetilde{S} / N$ if and only if $\gamma \in N$. More generally, if $\gamma N = \gamma_j^k N$ for some integer $k$, then $\gamma \cdot \widetilde{y}$ and $\widetilde{y}$ descend to the same component of $\widetilde{S} / N$.

With these preliminaries, we see from Proposition \ref{prop:quot} that if $\gamma \in \pi_1(S, x)$ then, for each $1 \leq j \leq t$, there is some integer $k$ so that $\gamma N = \gamma_j^k N$. As such, the boundary components of our quotient $\Sigma_\circ = \widetilde{S} / N$ are in one-to-one correspondence with those of $S$, as claimed in the lemma. Moreover, this shows that the action of the coset $\gamma_j N$ on the $j$th boundary circle $\Pi_\circ^{-1}(\partial D_j)$ of $\Sigma_\circ$ is by rotation by $2 \pi / n$ (since the action of $\gamma_j$ on each component of the universal cover is by translating by 1). We have taken our preferred generator of $\bb{Z}_n \cong \pi_1(S, x)/ N$ to be the coset $\gamma_1^{k_1} N$. Since $\gamma_1^{k_1} N = \gamma_j^{k_j}N$, this acts on the circle in $\Sigma_\circ$ over $\partial D_j$ by rotation by $e^{2 \pi \hat{i} k_j / n}$, as claimed.
\end{proof}

Fix $j$ and cap off the circle $\Pi_\circ^{-1}(\partial D_j) \subseteq \Sigma_\circ$ by gluing on a disk. Given our description of $\bb{Z}_n$-action on the boundary, this action extends across this disk to be the standard action of weight $k_j$. Thus, upon capping off this $j$th boundary component, we introduce a single fixed point with rotation number $k_j$. Doing this with all boundary components in $\Sigma_\circ$, we are left with a closed surface $\Sigma = \Sigma_\circ \cup \big(\sqcup_{j = 1}^t \bb{D}^2\big)$ that is equipped with an orientation-preserving homeomorphism 
\begin{equation}\label{eq:hdef}
h: \Sigma \longrightarrow \Sigma
\end{equation}
of order $n$ (given by the action of $\gamma_1^{k_1} N$) that has precisely the described fixed point data: the fixed points $x_1, \ldots, x_t$ are the centers of these disks we have capped, and the orbits of all other points have size $n$, so the associated $\bb{Z}_n$-action is indeed semifree. We can determine the genus by a simple Euler characteristic computation:
$$2 - 2g = \chi(\Sigma) = \chi(\Sigma_\circ) + t = n \chi(S) + t = n (\chi(S^2) - t) + t  = n(2 - t) + t.$$
This is equivalent to $g = (t - 2)(n-1)/2$. This gives the surface and $\bb{Z}_n$-action we were after.  

What remains is to determine the behavior of this action on the fundamental group of $\Sigma$. To obtain this, we pause to address some algebraic preliminaries about generators and relations of the group $N$ from above.

By Van Kampen's theorem, the effect of capping off the $j$th disk is to declare $\gamma_j^n$ to be the identity. That is, the fundamental group 
$$\pi_1(\Sigma, x_\circ) \cong N / M$$ 
is isomorphic to the quotient of $N$ by the normal subgroup $M$ in $N$ that is generated, as a normal subgroup, by $\gamma_j^n$ for $1 \leq j\leq t$. There are two issues with this description. The first is that we have been given generators for $N$, but only generators as a \emph{normal} subgroup, as opposed to generators as a group (and we prefer the latter when we pass to representation varieties). The second is that we will wish to express the relations $\gamma_j^n = e$ and $\gamma_1 \gamma_2 \ldots \gamma_t = e$ in terms of these generators we find. 

Now we turn to finding group generators for our subgroup $N$. Towards this end we set
$$ a_{j,\ell} \defeq \gamma_1^\ell \gamma_j^{-k_j} \gamma_1^{k_1} \gamma_1^{-\ell}  \hspace{.5cm} \textrm{and} \hspace{.5cm} d_{j, \ell} \defeq \gamma_1^\ell \gamma_j^n \gamma_1^{-\ell}$$
for $1 \leq j \leq t$ and $\ell \in \bb{Z}$. These are in $N$ since $N$ is normal, by definition. These are the generators we are after:

\begin{prop}\label{prop:Generators}
The group $N$ is generated, as a subgroup of $\pi_1(S)$, by the $a_{j, \ell}$ and $d_{j, \ell}$ for $1 \leq j \leq t$ and $\ell \in \bb{Z}$. The group $M$ is generated, as a subgroup of $\pi_1(S)$, by the $d_{j, \ell}$ for $1 \leq j \leq t$ and $\ell \in \bb{Z}$. 
\end{prop}

\begin{proof}
Write $N'$ for the group generated by the $a_{j, \ell}$ and $d_{j, \ell}$ for $1 \leq j \leq t$ and $\ell \in \bb{Z}$. We clearly have $N' \subseteq N$. For the reverse inclusion, it suffices to show that $N'$ is normal in $\pi_1(S)$. First, we plainly have $g N'  = N' g$ for $g = \gamma_1^\ell$ for any $\ell \in \bb{Z}$. Now if $g = \gamma_j^{k_j}$, then we have
$$\gamma_j^{k_j} N' = (a_{j, k_1})^{-1} \gamma_1^{k_1} N' = (a_{j, k_1})^{-1}  N' \gamma_1^{k_1} = N' (a_{j, k_1})^{-1}   \gamma_1^{k_1}  = N' \gamma_j^{k_j}.$$
Recalling that $k^{*}_j$ is the inverse of $k_j$ mod $n$, there is some integer $m_j$ so that $1 = k_j k^{*}_j + n m_j$, and so 
\begin{equation}\label{eq:gammajident}
\gamma_j  = (\gamma_j^{k_j})^{k^{*}_j} \gamma^{n m_j}.
\end{equation}
Then a proof similar to the one just given shows that
$$\gamma_j N' = N' \gamma_j$$
for all $j$. Thus, $g N' = N' g$ holds for all $g$ in the generating set of $\pi_1(S)$, and so holds for all $g \in \pi_1(S)$. This proves the claim about $N$, and the claim about $M$ is similar. 
\end{proof}

Of course, $a_{1, \ell} = e$ is trivial for all $\ell$ and so can be discarded from our generating set. Likewise, $d_{1, \ell} = d_{1, 0}$ and so we can discard all $d_{1, \ell}$ except $d_{1, 0}$. For $j \geq 2$, the $a_{j, \ell}$ and $d_{j, \ell}$ produce distinct words, though they do satisfy further nontrivial relations. For example, in $\pi_1(S)$, we have
$$\begin{array}{rcl}
\gamma_1^{n k_1} \gamma_j^{- n k_j} & = & (\gamma_1^{n k_1} \gamma_j^{-k_j} \gamma_1^{-(n-1) k_1} ) (\gamma_1^{(n-1)k_1} \gamma_j^{-k_j} \gamma_1^{-(n-2) k_1} )\times  \ldots\\
&& \ldots (\gamma_1^{2 k_1} \gamma_j^{-k_j} \gamma_1^{-k_1} ) (\gamma_1^{k_1} \gamma_j^{-k_j})\\
& = & a_{j, n k_1} a_{j, (n-1) k_1} \ldots a_{j, 2k_1} a_{j, k_1}
\end{array}$$
and so 
\begin{equation}\label{eq:rel1}
d_{1, 0}^{k_1} d_{j, 0}^{-k_j} = a_{j, n k_1} a_{j, (n-1) k_1} \ldots a_{j, 2k_1} a_{j, k_1} = \prod_{i = 0}^{n-1} a_{j, (n-i) k_1}.
\end{equation}

There is one more relation that we need to contend with, and this is coming from the relation $\gamma_1 \gamma_2 \ldots \gamma_t = e$ in $\pi_1(S)$. We need to express this in terms of the generators $a_{j, \ell}, d_{j, \ell}$, for $N$ and $M$, and this is ultimately supplied by Lemma \ref{lem:techlem0} below. Before stating the lemma, we begin with a sketch: The starting point is the identity (\ref{eq:gammajident}) from which we can express $\gamma_1 \gamma_2 \ldots \gamma_n$ first in terms of products of powers of the $\gamma_j^{k_j}$ and powers of $\gamma_j^n$. We are happy with the powers of $\gamma_j^n$, since they are generators of $M$. As for each power of $\gamma_j^{k_j}$, we can use 
$$\gamma_j^{k_j} = \gamma_1^{k_1}a_{j,0}^{-1}$$ 
to express this in terms of the generators $a_{j, 0}$ and various powers of $\gamma_1$. What remains are various powers of $\gamma_1^{k_1}$, separated by $a_{j, \ell}$ terms and $d_{j, \ell}$ terms. Since $N$ is normal in $\pi_1(S)$, we can consolidate these copies of $\gamma_1^{k_1}$ using
\begin{equation}\label{eq:normalc}
a^{-1}_{j, \ell} \gamma^i_1= \gamma_1^i a_{j, \ell -i}^{-1} 
\end{equation}
to get a single copy of $\gamma_1$ to some power followed by a big product of our desired generators $a_{j, \ell}$ and $d_{j, \ell}$. Tracking the exponents, this power of $\gamma_1$ will be a multiple of $\sum_j k^{*}_j \in n \bb{Z}$, and so is a power of $\gamma_1^n = d_{1, 0}$, which is in $M$.

To state our lemma, let $Fr(\gamma_1, \ldots, \gamma_t)$ be the free group on the generators $\gamma_1, \ldots, \gamma_t$. In the statement of lemma, the symbols $a_{j, \ell}$ and $d_{j, \ell}$ should be defined as above, but viewed as elements of $Fr(\gamma_1, \ldots, \gamma_t)$.

\begin{lemma}\label{lem:techlem0}
For each $\mu \in \bb{Z}$, there is a set 
$$\left\{(j(\mu; k), l(\mu; k)) \right\}_{k = 1}^{K_\mu}\subseteq \left\{2, \ldots, t \right\} \times \bb{Z}$$ 
of ordered pairs of integers so that
\begin{equation}\label{eq:rel2}
\prod_{j = 2}^t \prod_{i = 1}^{k^{*}_j} a^{-1}_{j, \ell(i, j)+ \mu}  = \gamma_1^\mu (\gamma_1 \gamma_2 \ldots \gamma_t) \gamma_1^{-\mu} \prod_{k = 1}^{K_\mu} d_{j(\mu; k), l(\mu; k)},
\end{equation}
in $Fr(\gamma_1, \ldots, \gamma_t)$, where the $\ell(i, j)$ are defined in (\ref{eq:elldef}).
\end{lemma}

Ultimately we will pass to the $M$-cosets, which will kill off all factors $d_{j, \ell}$; as such, the precise values of $j(\mu; k)$ and $l(\mu; k)$ will not be relevant to us.

\begin{proof}[Proof of Lemma \ref{lem:techlem0}]
Let $M'$ be the subgroup of $Fr(\gamma_1, \ldots, \gamma_t)$ generated by the $d_{j, \ell}$ (this descends to our normal subgroup $M$ under the obvious surjection $Fr(\gamma_1, \ldots, \gamma_t) \rightarrow \pi_1(S, x)$). The argument of Proposition \ref{prop:Generators} carries over to this situation to show that $M'$ is normal in $Fr(\gamma_1, \ldots, \gamma_t)$. 

Working with $M'$-cosets, we have
$$\begin{array}{rcl}
(\gamma_1 \gamma_2 \ldots \gamma_t)M' & = & \gamma_1 (\gamma_2^{k_2})^{k^{*}_2} \ldots (\gamma_t^{k_t})^{k^{*}_t} M'\\
& = & \gamma_1 (\gamma_1^{k_1} a_{2, 0}^{-1})^{k^{*}_2}(\gamma_1^{k_1} a_{3, 0}^{-1})^{k^{*}_3} \ldots (\gamma_1^{k_1} a_{t, 0}^{-1})^{k^{*}_t} M'.
\end{array}$$
To evaluate this further, use (\ref{eq:normalc}) to get
$$(\gamma_1^{k_1} a_{j, 0}^{-1})^{k^{*}_j} = \gamma_1^{k^{*}_j k_1} a^{-1}_{j, -(k^{*}_j -1) k_1} a^{-1}_{j, -(k^{*}_j -2) k_1}  \ldots a^{-1}_{j, 0}.$$
Thus
$$\begin{array}{rcl}
(\gamma_1 \gamma_2 \ldots \gamma_t)M' & = & \gamma_1 \gamma_1^{k^{*}_2 k_1} \Big( a^{-1}_{2, -(k^{*}_2 -1) k_1} a^{-1}_{2, -(k^{*}_2-2) k_1}  \ldots a^{-1}_{2, 0} \Big)\\
&& \times \gamma_1^{k^{*}_3 k_1} \Big( a^{-1}_{3, -(k^{*}_3 -1) k_1} a^{-1}_{3, -(k^{*}_3-2) k_1}  \ldots a^{-1}_{3, 0} \Big)\\
&& \times  \ldots \\
&& \times \gamma_1^{k^{*}_t k_1} \Big( a^{-1}_{t, -(k^{*}_t -1) k_1} a^{-1}_{t, -(k^{*}_t-2) k_1}  \ldots a^{-1}_{t, 0} \Big)M'.
\end{array}$$
Now use (\ref{eq:normalc}) again to move the copies of $\gamma_1$ all to the front to get
$$\begin{array}{rcl}
(\gamma_1 \gamma_2 \ldots \gamma_t)M'  & = & \Big( \gamma_1^{1 + k^{*}_2 k_1 + k^{*}_3 k_1 + \ldots + k^{*}_t k_1}  \prod_{j = 2}^t \prod_{i = 1}^{k^{*}_j} a^{-1}_{j, \ell(i, j)}\Big) M'
\end{array}$$
where $\ell(i, j)$ is as in (\ref{eq:elldef}). Recall our consistency condition on the $k^{*}_j$ gives
$$\begin{array}{rcl}
1 + k^{*}_2 k_1 + k^{*}_3 k_1 + \ldots + k^{*}_t k_1 & \equiv & k^{*}_1 k_1 + k^{*}_2 k_1 + k^{*}_3 k_1 + \ldots + k^{*}_t k_1 \\
&= & k_1 \sum_{j} k^{*}_j \\
& \equiv & 0
\end{array}$$
where the equivalence is mod $n$. Thus, $\gamma_1^{1 + k^{*}_2 k_1 + k^{*}_3 k_1 + \ldots + k^{*}_t k_1}  \in M'$ and so we conclude
$$ \prod_{j = 2}^t \prod_{i = 1}^{k^{*}_j} a^{-1}_{j, \ell(i, j)}  M' = (\gamma_1 \gamma_2 \ldots \gamma_t) M'.$$
Since $a_{j, \ell + \mu} = \gamma_1^\mu a_{j, \ell} \gamma_1^{-\mu}$, and $M'$ is normal, this implies
$$\Big( \prod_{j = 2}^t \prod_{i = 1}^{k^{*}_j} a^{-1}_{j, \ell(i, j)+\mu} \Big) \gamma_1^\mu (\gamma_1 \gamma_2 \ldots \gamma_t)^{-1} \gamma_1^{-\mu} \in  M'$$
for all $\mu \in \bb{Z}$. The group $M'$ is generated by the $d_{j, \ell}$, so we can write this term on the left as a product of the $d_{j, \ell}$, as desired.
\end{proof}

By Proposition \ref{prop:Generators}, we see that the cosets $a_{j, \ell} M$ generate $\pi_1(\Sigma, x_\circ) \cong N/ M$, and (\ref{eq:rel1}) and (\ref{eq:rel2}) implies they satisfy 
\begin{equation}\label{eq:rel12}
\Bigg(\prod_{i = 0}^{n-1} a_{j, (n-i) k_1} \Bigg) M= M, \hspace{.5cm} \textrm{and} \hspace{.5cm} \Bigg( \prod_{j = 2}^t \prod_{i = 1}^{k^{*}_j} a^{-1}_{j, \ell(i, j)+ \mu}\Bigg) M = M
\end{equation}
for all $\mu \in \bb{Z}$. Since $\gamma_1^n \in M$ and $M$ is normal, we have
$$a_{j, \ell + n} M = a_{j, \ell} M$$
are equal for all $\ell \in \bb{Z}$ and $1 \leq j \leq t$. As such, when considering the $M$-cosets $a_{j, \ell} M$, we can treat the index $\ell \in \bb{Z}_n$ mod $n$. 

Our remaining task is to show that (\ref{eq:rel12}) generate all relations satisfied by the $a_{j, \ell} M$ for $2 \leq j \leq t$ and $\ell \in \bb{Z}_n$. To do this, let $Fr(\alpha_{j, \ell})$ be the free group generated by new symbols $\alpha_{j, \ell}$ indexed by $2 \leq j \leq t $ and $\ell \in \bb{Z}$. Consider the map
\begin{equation}\label{eq:Psi}
\Psi: Fr(\alpha_{j, \ell} ) \longrightarrow N/M
\end{equation}
given by sending $\alpha_{j, \ell}$ to $a_{j, \ell}^{-1} M$. In Section \ref{sec:Completing}, we will show that the kernel of this map is generated by the relations claimed in Theorem \ref{thm:Action}, which are those of (\ref{eq:rel12}) expressed in terms of the $\alpha_{j, \ell}$. To do this, we find it convenient to bring the $d_{j, \ell}$ back into the picture, as we will need to work with some explicit formulas. We break the task up into two steps, first treating the relations from (\ref{eq:rel12}) associated to (\ref{eq:rel1}), and then those associated to (\ref{eq:rel2}). 

\subsection{Step 1: Addressing the relations associated to (\ref{eq:rel1})}

Let $Fr(a_{j, \ell}, d_{j, \ell})$ be the free group generated by $d_{1, 0}$ and also the $a_{j, \ell}$ and $d_{j, \ell}$ for $2 \leq j \leq t$ and $\ell \in \bb{Z}$. Throughout this step, we will view the $a_{j, \ell}$ and $d_{j, \ell}$ as being elements of $Fr(a_{j, \ell}, d_{j, \ell})$ as opposed to $\pi_1(S, x)$; thus, here, the $a_{j, \ell}, d_{j, \ell}$ do not satisfy any nontrivial relations. Likewise, throughout this step will view the $\gamma_j$ as being elements of the free group $Fr(\gamma_1, \ldots, \gamma_t)$; it won't be until after this step that we return to viewing the $\gamma_j$ as elements of $\pi_1(S) \cong Fr(\gamma_1, \ldots, \gamma_t) / \langle \gamma_1 \gamma_2 \ldots \gamma_t \rangle_{norm}$, where $\langle \gamma_1 \gamma_2 \ldots \gamma_t \rangle_{norm}$ is the normal subgroup generated by the product $\gamma_1 \gamma_2 \ldots \gamma_t$. Consider the homomorphism
$$\psi: Fr(a_{j, \ell}, d_{j, \ell}) \longrightarrow Fr(\gamma_1, \ldots, \gamma_t)$$
sending $a_{j, \ell} \mapsto \gamma_1^{\ell} \gamma_j^{-k_j} \gamma_1^{k_1 - \ell}$ and $d_{j, \ell} \mapsto \gamma_1^{\ell} \gamma_j^n \gamma_1^{-\ell}$. Here is the main technical result of this step.

\begin{lemma}\label{lem:techlem1}
The kernel of $\psi$ is generated, as a normal subgroup, by the words 
$$W_{j, r} \defeq d_{1, 0}^{r k_1} d_{j, 0}^{-r k_j} \prod_{i = 0}^{r n -1} a^{-1}_{j, i k_1}$$
$$A_{j, \ell} \defeq d_{1, 0} a_{j, \ell} d_{1, 0}^{-1} a_{j, \ell+n}^{-1}, \hspace{.5cm}  D_{j, \ell} \defeq d_{1, 0} d_{j, \ell} d_{1, 0}^{-1} d_{j, \ell+n}^{-1}$$
for all $j \in \left\{2, \ldots, t \right\}$, all $r \in \bb{Z}$, and all $\ell \in \bb{Z}$. 
\end{lemma}

\begin{proof}
The proof strategy of (\ref{eq:rel1}) carries over to show that the words $W_{j, r}$ are all contained in the kernel of $\psi$ (i.e., the proof never used the relation $\gamma_1 \gamma_2 \ldots \gamma_t = e$ that holds in $\pi_1(S)$, but is not valid here). Since 
$$\begin{array}{rcl}
\psi(d_{1, 0} a_{j, \ell}) & = & \gamma_1^n \gamma^\ell \gamma_j^{-k_j} \gamma_1^{k_1 - \ell}\\
& = & \gamma^{\ell + n} \gamma_j^{-k_j} \gamma_1^{k_1 - (\ell+n)} \gamma_1^n\\
& = & \psi(a_{j, \ell+n} d_{1, 0})
\end{array}$$
it follows that $A_{j, \ell} \in \ker(\psi)$. A similar computation shows that $D_{j, \ell}  \in \ker(\psi)$. 

For the reverse inclusion, suppose $w \in \ker(\psi)$ is a nontrivial element of this kernel. Write $w$ as a reduced word in the generators $d_{1, 0}, a_{j, \ell}$, and $d_{j, \ell}$. First, $w$ cannot consist entirely of terms of the form $a_{j, \ell}$ for various $j, \ell$, since that would induce a nontrivial word in the $\gamma_j$, and there is no such word. Likewise, $w$ cannot consist entirely of the $d_{j, \ell}$ for the same reason. By this same type of argument, there is some $j_0 \in \left\{2, \ldots, t \right\}$ and $\ell, \ell' \in \bb{Z}$ so that $a_{j_0, \ell}$ and $d_{j_0, \ell'}$ are terms appearing in $w$. For this $j_0$, enumerate all of these occurrences as $x_{j_0, \ell_1} , x_{j_0, \ell_2} , \ldots , x_{j_0, \ell_n}$ where $x_{j_0, \ell_k} \in \left\{ a_{j_0, \ell_k}, a^{-1}_{j_0, \ell_k}, d_{j_0, \ell_k}, d^{-1}_{j_0, \ell_k} \right\}$ and the ordering corresponds to the order in which these appear in our reduced word $w$. Then
$$e = \psi(w) = \ldots \psi (x_{j_0, \ell_1} ) \ldots \psi(x_{j_0, \ell_2}) \ldots \psi(x_{j_0, \ell_n}) \ldots.$$
The only appearances of $\gamma_{j_0}$ are those in $\psi (x_{j_0, \ell_k} )$, and they are all flanked by powers of $\gamma_1$ (some of these powers may be zero). Since the $\gamma_j$ are free in the codomain of $\psi$, these terms $\gamma(x_{j_0, \ell_k})$ are either directly adjacent to each other or are separated only by a power of $\gamma_1 = \psi(d_{1, 0})$. That is, we can write
$$w = \ldots (x_{j_0, \ell_1} d_{1, 0}^{n_1} x_{j_0, \ell_2} d_{1, 0}^{n_2}  \ldots x_{j_0, \ell_n}d_{1, 0}^{n_{n}}) \ldots$$
for some powers $n_k \in \bb{Z}$. By repeated use of 
$$\psi(d_{1, 0} a_{j_0, \ell_k}) = \psi(a_{j_0, \ell_k + n} d_{1, 0}), \hspace{.5cm} \psi(d_{1, 0} d_{j_0, \ell_k}) = \psi(d_{j_0, \ell_k + n} d_{1, 0}),$$
we can reduce to the case where all powers $n_1 = n_2 = \ldots n_n = 0$ vanish, at the cost of having some additional power of $d_{1, 0}$ on the far right side. In this way, we have reduced to the case where $w$ contains the subword
$$w' \defeq x_{j_0, \ell_1}  x_{j_0, \ell_2}  \ldots x_{j_0, \ell_n}$$
where each $x$ is one of $a, a^{-1}, d$, or $d^{-1}$. By construction, this subword $w'$ contains all occurrences in $w$ of the symbols $a_{j_0, \ell}, a^{-1}_{j_0, \ell},d_{j_0, \ell}, d_{j_0, \ell}^{-1}$ for any $\ell$. In $Fr(\gamma_1, \ldots, \gamma_t)$ if $j$ and $j'$ are distinct in $\left\{2, \ldots, t \right\}$, then there are no nontrivial relations among the $\gamma_1^{\ell} \gamma_{j}^{-k_j} \gamma_1^{k_1 - \ell}$ and $\gamma_1^{\ell} \gamma_{j'}^n \gamma_1^{-\ell}$ for any $\ell$. It follows that, when plugged into $w'$, all terms cancel except possibly some residual copies of $\gamma_1$. That is, there is some integer $h$ so that 
$$\psi(w'd_{1, 0}^h) = e.$$
It suffices to show that $w' d_{1, 0}^h$ is in the normal subgroup generated by the $A_{j, \ell}, D_{j, \ell}$, and $W_{j, r}$. 

As with $w$ above, we can assume that $w' = x_{j_0, \ell_1}  x_{j_0, \ell_2}  \ldots x_{j_0, \ell_n}$ is in reduced form (though $w$ was in reduced form, this may have been destroyed for $w'$ through the reductions above). If this makes $w' d_{1, 0}^h$ trivial, then we are done, so we can assume that $w' d_{1, 0}^h$ is a nontrivial word; thus $w'$ is a nontrivial word since $d_{1, 0}$ is not in the kernel of $\psi$. By conjugating $w' d_{1, 0}^h$ by $d_{j_0, \ell}$ terms for various $\ell$, and using $\psi(d_{1, 0} d_{j_0, \ell}) = \psi(d_{j_0, \ell+n} d_{1, 0})$, we may assume that the first term $x_{j_0, \ell_1} = a^{\pm 1}_{j_0, \ell_1}$ in $w'$ is an $a$-term. We will work in the case where $x_{j_0, \ell_1} = a_{j_0, \ell_1}$; the other case is similar (it can be reduced to this case by a combination of taking the inverse of $w'$ and a judicious choice of conjugation). 

Let $f$ be the largest positive integer so that $a_{j_0, \ell_1} a^{\pm 1}_{j_0, \ell_2} \ldots a^{\pm 1}_{j_0, \ell_f}$ is the starting sequence of $w'$. As we have seen with $w$ above, $w'$ cannot consist of $a$-terms only, so by the maximality of $f$ we can write
$$w' = a_{j_0, \ell_1} a^{\pm 1}_{j_0, \ell_2} \ldots a^{\pm 1}_{j_0, \ell_f} d^{\pm 1}_{j_0, \ell_{f+1}} \ldots.$$
Zooming in on the first two terms, $a_{j_0, \ell_1} a^{\pm 1}_{j_0, \ell_2}$, assume that the $\pm$ here is $-$. Then we would have 
$$\psi(a_{j_0, \ell_1} a^{\pm 1}_{j_0, \ell_2} )= \psi(a_{j_0, \ell_1} a^{-1}_{j_0, \ell_2}) = \gamma_1^{k_1} \gamma_{j_0}^{-k_{j_0}} \gamma_1^{\ell_2 - \ell_1} \gamma_{j_0}^{k_{j_0}} \gamma_1^{- \ell_2}.$$
Since $j_0 \neq 1$ and the $\gamma_j$ are free, this forces $\ell_1 = \ell_2$ otherwise there is no way for $w' d^h_{1, 0}$ to map to the identity under $\psi$. However, $w'$ is not in reduced form when $\ell_1 = \ell_2$. We must therefore have that $\pm = +$ in this first term, from which we similarly compute
$$\psi(a_{j_0, \ell_1} a^{\pm 1}_{j_0, \ell_2} ) = \psi(a_{j_0, \ell_1} a^{+ 1}_{j_0, \ell_2} ) = \gamma_1^{\ell_1} \gamma_{j_0}^{-k_{j_0}} \gamma_1^{k_1 - \ell_1 + \ell_2} \gamma_{j_0}^{-k_{j_0}} \gamma_{1}^{k_1 - \ell_2},$$
and deduce that $\ell_2 = \ell_1 - k_1$. The same observation holds for the other $a$-terms in this first sequence of $w'$: The $\pm 1$'s are all $+1$'s and
$$\ell_{k+1} = \ell_k - k_1$$
for all $1 \leq k \leq f-1$; that is, $\ell_{k+1} = \ell_1 - k k_1$ for $0 \leq k \leq f-1$. By conjugating $w' d_{1, 0}^h$ again, we can make the $\ell_1$ be anything we want at the cost of changing the exponent $h$ (which we can tolerate). We find it convenient to do this so that $\ell_1 = f k_1$. Thus we have reduced to the case where $w'$ has the form 
$$w' = a_{j_0, f k_1} a_{j_0, (f - 1) k_1} \ldots a_{j_0, k_1 } d^{\pm 1}_{j_0, \ell_{f+1}} \ldots.$$

We want to next take the maximal string of $d$-terms starting at $d^{\pm 1}_{j_0, 0}$. If there is no next term, or the next term is an $a$-term, then this maximal string consists just of $d^{\pm 1}_{j_0, \ell_{f+1}}$. Suppose that the next term is a $d$-term, so
$$w' = a_{j_0, f k_1} a_{j_0, (f-1) k_1} \ldots a_{j_0, k_1 } d^{\pm 1}_{j_0, \ell_{t+1}} d^{\pm 1}_{j_0, \ell_{f+2}} \ldots.$$
Plugging $a_{j_0, k_1 } d^{\pm 1}_{j_0, \ell_{f+1}} d^{\pm 1}_{j_0, \ell_{f+2}}$ into $\psi$, we see that $\ell_{f+1} = \ell_{f+2} = 0$. A similar computation holds for all adjacent $d$-terms. We deduce that, in general, 
$$w' = a_{j_0, f k_1} a_{j_0, (f-1) k_1} \ldots a_{j_0, k_1 } d^{m}_{j_0, 0} a_{j_0, 0}^{m'}\ldots$$
for some integers $m, m', s$ with $m$ non-zero (and $m' = 0$ if and only if $w'$ has no terms after $d^{m}_{j_0, 0} $). Repeating down the line, there are integers $f_1, \ldots, f_T \geq 1$ and non-zero integers $m_1, \ldots, m_T$ so that $w' d_{1, 0}^h$ equals
$$ \begin{array}{l}
(a_{j_0, f_1 k_1} \ldots a_{j_0, k_1 } d^{m_1}_{j_0, 0} ) (a_{j_0, f_2 k_1} \ldots a_{j_0, k_1 } d^{m_2}_{j_0, 0} ) \times \\
\; \hspace{4cm} \ldots \times (a_{j_0, e_T k_1} \ldots a_{j_0, k_1 } d^{m_T}_{j_0, 0} ) d_{1, 0}^{h'}\end{array}$$
for some $h' \in \bb{Z}$. The lemma follows from the next claim.

\medskip

\noindent \emph{Claim: For all $1 \leq i \leq T$ there is some integer $r_i$ so that $f_i = r_i n$, $m_i = r_i k_{j_0}$, and $h' = -k_1 (r_1 + \ldots + r_T)$. Thus
$$w' d_{1, 0}^h = \Big(W_{j_0, r_1}^{-1} d_{1, 0}^{k_1 r_1} \Big) \ldots \Big( W_{j_0, r_T}^{-1} d_{1, 0}^{k_1 r_T} \Big) d_{1, 0}^{- k_1(r_1 + \ldots + r_T)} \in \ker(\psi).$$}

\medskip

When $T = 1$ we have
$$e = \psi(w' d_{1, 0}^h) = \gamma_1^{f_1 k_1} \gamma_{j_0}^{-k_{j_0} f_1 + n m_1} \gamma_1^{n h'}$$
from which we deduce that $n m_1 = k_{j_0} f_1$ and $n h' = - f_1 k_1$. Since $k_{j_0}$ is invertible mod $n$, the equation $n m_1 = k_{j_0} f_1$ implies that $f_1 = r_1 n$ is a multiple of $n$ and thus $m_1 = r_1 k_{j_0} $ is the same multiple of $k_{j_0}$. Then $n h' = - f_1 k_1$ implies $h' = - k_1 r_1$, so this proves the result for $T = 1$. 

Working inductively, assume the claim holds for $T \geq 1$. Then for each $1 \leq i \leq T$ we have
$$a_{j_0, f_i k_1} \ldots a_{j_0, k_1 } d^{m_i}_{j_0, 0}  = W_{j_0, r_i}^{-1} d_{1, 0}^{k_1 r_i}.$$
Since $\psi(w' d_{1, 0}^h) = e$ is the identity, this gives
$$e =  \psi(W_{j_0, r_1}^{-1} d_{1, 0}^{k_1 r_1}) \ldots \psi(W_{j_0, r_T}^{-1} d_{1, 0}^{k_1 r_T}) \gamma_1^{f_{T+1} k_1} \gamma_{j_0}^{-k_{j_0} f_{T+1} + n m_{T+1}} \gamma_1^{n h'}.$$
Since we know the $W_{j, r}$ are in the kernel of $\psi$, we can write this as 
$$e = \gamma_{1}^{n k_1(r_1 + \ldots + r_T) } \gamma_1^{f_{T+1} k_1} \gamma_{j_0}^{-k_{j_0} f_{T+1} + n m_{T+1}} \gamma_1^{n h'}.$$
The same argument as in the base case shows that there is some integer $r_{T+1}$ so that $f_{T+1} = r_{T+1} n$ and $m_{T+1} = r_{T+1} k_{j_0}$, and $h' = -k_1(r_1 + \ldots + r_{T+1})$. 
\end{proof}

\subsection{Step 2: Addressing the relations associated to (\ref{eq:rel2})}

We will continue to use the notation from the previous section. Here we will be interested in the inverse image $\psi^{-1}(\langle \gamma_1 \gamma_2 \ldots \gamma_t \rangle_{norm} )$ of the normal subgroup of $Fr(\gamma_1, \ldots, \gamma_t)$ generated by the product $\gamma_1 \gamma_2 \ldots \gamma_t$.

\begin{lemma}\label{lem:techlem2}
The inverse image $\psi^{-1}(\langle \gamma_1 \gamma_2 \ldots \gamma_t \rangle_{norm} )$ is generated, as a normal subgroup, by the words $W_{j, r}, A_{j, \ell}, D_{j, \ell}$, and the additional words
$$P_{\mu} \defeq \Big(\prod_{j = 2}^t \prod_{i = 1}^{k^{*}_j} a^{-1}_{j, \ell(i, j) + \mu} \Big)\Big(\prod_{k= 1}^{K_\mu} d_{j(\mu; k), \ell(\mu; k)}\Big)^{-1}$$
for $2 \leq j \leq t$ and $r , \ell, \mu \in \bb{Z}$.
\end{lemma}

\begin{proof}
By Lemmas \ref{lem:techlem1}, and \ref{lem:techlem0} the words $W_{j, r}, A_{j, \ell}, D_{j, \ell}$, and $P_\mu$ are indeed contained in this inverse image. To see the reverse direction, it suffices to show that if $w \in Fr(a_{j, \ell}, d_{j, \ell})$ is a word with $\psi(w)$ conjugate to $\epsilon \defeq \gamma_1 \gamma_2 \ldots \gamma_t$, then there is some $\mu \in \bb{Z}$ and $u \in Fr(a_{j, \ell}, d_{j, \ell})$ so that the cosets
$$w K = u P_\mu u^{-1} K$$
are equal, where $K \defeq \ker(\psi)$. 

As a first case, assume $\psi(w)= \gamma_1^{\mu} \epsilon \gamma_1^{-\mu}$ for some $\mu \in \bb{Z}$. Since $\psi(P_\mu) = \gamma_1^{\mu} \epsilon \gamma_1^{-\mu}$, we have $w K = P_\mu K$ as desired. 

For the next case, assume $\psi(w) = \psi(u) \epsilon \psi(u)^{-1}$ for some $u \in Fr(a_{j, \ell}, d_{j, \ell})$. Then $\psi(u P_0 u^{-1}) = \psi(w)$ and so $wK = uP_0 u^{-1} K$.

As a third case, suppose $\psi(w) = \gamma_j \epsilon \gamma_j^{-1}$ for some $2\leq j \leq t$. We have
$$\gamma_j = (\gamma_j^{k_j})^{k^{*}_j} \gamma_j^{n_j n} = (\gamma_1^{k_1} \psi(a_{j, 0}^{-1}))^{k^{*}_j} \psi(d_{j, 0})^{n_j}.$$
This is a product of terms from the previous two cases, so this case reduces to those.

Finally, if $\psi(w) = v \epsilon v^{-1}$ for an arbitrary $v \in Fr(\gamma_1, \ldots, \gamma_t)$, then by expressing $v$ as a product of the $\gamma_j$, we are reduced to the first and third cases.

\end{proof}

\subsection{The generators and relations of Theorem \ref{thm:Action}}\label{sec:Completing}

In an attempt to avoid confusion in reference to the steps above, for the rest of the proof we will reserve the symbols $a_{j, \ell}$ and $d_{j, \ell}$ for the generators of the free group $Fr(a_{j, \ell}, d_{j, \ell})$, reserving the more cumbersome notation $\gamma_1^\ell \gamma_{j}^{-k_j} \gamma_1^{k_1 - \ell}$ and $\gamma^\ell_1 \gamma_j^n \gamma_1^{-\ell}$ for their images under $\psi$ in $Fr(\gamma_1, \ldots, \gamma_t)$, and then the cosets 
$$\gamma_1^\ell \gamma_{j}^{-k_j} \gamma_1^{k_1 - \ell} \langle \gamma_1 \gamma_2 \ldots \gamma_t \rangle_{norm}\hspace{.5cm} \textrm{and} \hspace{.5cm} \gamma^\ell_1 \gamma_j^n \gamma_1^{-\ell} \langle \gamma_1 \gamma_2 \ldots \gamma_t \rangle_{norm}$$ 
for the images of these in $\pi_1(S, x) \cong Fr(\gamma_1, \ldots, \gamma_t) / \langle \gamma_1 \gamma_2 \ldots \gamma_t \rangle_{norm}$. 

Consider the composition
$$Fr(a_{j, \ell}, d_{j, \ell}) \stackrel{\psi}{\longrightarrow} Fr(\gamma_1, \dots, \gamma_t) \longrightarrow \frac{Fr(\gamma_1, \dots, \gamma_t)}{\langle \gamma_1 \gamma_2 \ldots \gamma_t \rangle_{norm}} \cong \pi_1(S, x).$$
By Proposition \ref{prop:Generators}, this has image the subgroup $N \leq \pi_1(S, x)$. Write
$$\phi: Fr(a_{j, \ell}, d_{j, \ell})  \longrightarrow N$$
for the map obtained from $\psi$ by corestriction. From Lemma \ref{lem:techlem2} we see that the kernel of $\phi$ is generated by the $W_{j, r}, A_{j,\ell}, D_{j, \ell}$, and $P_\mu$. Let $M'$ be the subgroup of $Fr(a_{j, \ell}, d_{j, \ell})$ normally generated by $d_{1, 0}$ and the $d_{j, \ell}$ for $2 \leq j \leq t$ and $\ell \in \bb{Z}$. Essentially by definition of the maps involved, the image of $M'$ under $\phi$ is precisely the normal subgroup $M$. As such, $\phi$ descends to a surjective homomorphism
$$\Phi: Fr(a_{j, \ell}, d_{j, \ell})/ M' \longrightarrow N / M.$$
We thus have an isomorphism
$$\frac{Fr(a_{j, \ell}, d_{j, \ell})/M'}{\ker(\Phi)} \cong N/ M.$$
By definition of $M'$, the quotient $Fr(a_{j, \ell}, d_{j, \ell}) / M'$ is freely generated by the cosets $a_{j, \ell} M'$ for $2 \leq j \leq t$ and $\ell \in \bb{Z}$. From our understanding of the kernel of $\phi$, we see that the kernel of $\Phi$ is generated by the cosets 
$$\begin{array}{rcl}
W_{j, r} M' & = & \prod_{i = 0}^{r n - 1} a_{j, i k_1}^{-1} M'\\
A_{j, \ell} M' & = & a_{j, \ell} a_{j, \ell+n}^{-1} M'\\
D_{j, \ell} M' &= & M'\\
P_\mu M' & = & \prod_{j = 2}^t \prod_{i = 1}^{k^{*}_j} a^{-1}_{j, \ell(i, j) + \mu}M'.
\end{array}$$
A priori we need to keep track of these for $2 \leq j \leq t$ and $r, \ell, \mu \in \bb{Z}$, but we will see that there are many redundancies and simplifications that can be made. To start, we can plainly drop the relations $D_{j, \ell} M'$ since these equal the trivial element $M'$ in $Fr(a_{j, \ell}, d_{j, \ell})/ M'$. The relations $a_{j, \ell} a_{j, \ell+n}^{-1} M'$ give our biggest reduction, since they show the index $\ell$ can be treated mod $n$. Indeed, setting $\alpha_{j, \ell}$ equal to the $\ker(\Phi)$-coset of $a_{j, \ell}^{-1} M'$, we have
$$\alpha_{j, \ell} = \alpha_{j, \ell + n}$$
and these generate $\frac{Fr(a_{j, \ell}, d_{j, \ell})/M'}{\ker(\Phi)}$ for $2 \leq j \leq t$ and $\ell \in \bb{Z}_n$. In terms of these generators, the remaining relations become
$$\prod_{i = 0}^{r n - 1} \alpha_{j, i k_1}, \hspace{.5cm} \prod_{j = 2}^t \prod_{i = 1}^{k^{*}_j} \alpha_{j, \ell(i, j) + \mu}.$$
By the mod $n$ periodicity in the second index of the $\alpha_{j, \ell}$, we need only consider these for $r = 1$ and $\mu \in \bb{Z}_n$, as all other values of $r, \mu \in \bb{Z}$ reduce to these. Thus $\pi_1(\Sigma, x_\circ)$ has the presentation as claimed in Theorem \ref{thm:Action}, so the same is true of $\pi_1(\Sigma, x_1) \cong \pi_1(\Sigma, x_\circ)$. What remains is to understand the action of $\bb{Z}_n$ on these generators.

\subsection{The action on \texorpdfstring{$\pi_1(\Sigma)$}{}}

First we give a description of the action from a purely algebraic perspective. Then we consider the actual periodic map $h: \Sigma \rightarrow \Sigma$ from (\ref{eq:hdef}) and show that its pushforward induces this algebraic action. 

In general, suppose $\pi$ is a group and $H$ is a normal subgroup, and consider the short exact sequence
$$e \longrightarrow H \longrightarrow \pi \longrightarrow \pi / H \longrightarrow e.$$
Note that any right splitting (if it exists) of this sequence induces a right action of $\pi / H$ on $H$. Indeed, if $\sigma: \pi/H \rightarrow \pi$ is such a right splitting, then this action is defined by
$$z \cdot h \defeq \sigma(z)^{-1} h \sigma(z)$$
for $z \in \pi/ H$ and $h \in H$, where the product on the right is taking place in $\pi$. 

Turning now to our situation, we start with the short exact sequence
$$e \longrightarrow N \longrightarrow \pi_1(S, x) \longrightarrow \pi_1(S, x) / N \longrightarrow e.$$
Since $M$ is a normal subgroup of $N$, we can quotient all terms by the natural $M$-action to get a second short exact sequence
\begin{equation}\label{eq:ses}
e \longrightarrow N/M \longrightarrow \pi_1(S, x)/M \longrightarrow \pi_1(S, x) / N \longrightarrow e.
\end{equation}
Recall from (\ref{eq:PrefIso}) that our preferred generator of $\pi_1(S, x) / N \cong \bb{Z}_n$ is $\gamma_1^{k_1} N$. Since $\gamma_1^{k_1} M$ has order $n$ in $\pi_1(S, x)/M$, the assignment
$$\sigma(\gamma_1^{k_1} N) \defeq \gamma_1^{k_1} M$$
extends to a well-defined splitting of (\ref{eq:ses}). It is in this way that $\bb{Z}_n$ acts on $N/M$---by conjugation by powers of $\gamma_1^{k_1} M$. The action of the generator $\gamma_1^{k_1} N$ of $\pi_1(S, x) / N$ is the map on $\pi_1(\Sigma, x_\circ) \cong \pi_1(S, x) / M$ sending
\begin{equation}\label{eq:actionform}
a_{j, \ell}M \longmapsto a_{j, \ell - k_1}M
\end{equation}
for all $j$ and all $\ell$, which structurally recovers the action as claimed in Theorem \ref{thm:Action} (though notice that the basepoint is off; we will deal with this momentarily). 

To finish the proof, we will show that the period-$n$ homeomorphism $h: \Sigma \rightarrow \Sigma$ from (\ref{eq:hdef}) induces precisely this map on the fundamental group. There is a subtlety here with which we must contend: The various groups considered above are all defined relative to the basepoint $x_\circ \in \Sigma_\circ \subseteq \Sigma$, but this is \emph{not} fixed by $h$. In the next paragraph we will show that the composition of the map
$$h_*: \pi_1(\Sigma, x_\circ) \longrightarrow \pi_1(\Sigma, h(x_\circ)),$$
with the change-of-basepoint map $\pi_1(\Sigma, h(x_\circ)) \cong \pi_1(\Sigma, x_\circ)$ has the form (\ref{eq:actionform}). Then the claimed formula for $h_*: \pi_1(\Sigma, x_1) \rightarrow \pi_1(\Sigma, x_1)$ will follow again from another use of the change-of-basepoint isomorphism (this time from $x_\circ$ to $x_1$).

Towards this end, let $\eta \in \pi_1(\Sigma, x_\circ)$; we will not distinguish in the notation between $\eta$ and a loop $[0, 1] \rightarrow \Sigma$ representing this homotopy class. By homotoping $\eta$, we can arrange so that it is a loop in the subspace $\Sigma_\circ = \widetilde{S}/ N \subseteq \Sigma$. Let $\widetilde{\eta}: [0, 1] \rightarrow \widetilde{S}$ be a lift of $\eta$ to the universal cover $\widetilde{S}$. By definition, $\phi \circ \eta$ is the image in $\widetilde{S}/ N \subseteq \Sigma$ of the path $\gamma_1^{k_1} \cdot \widetilde{\eta}$, where $\gamma_1^{k_1} \in \pi_1(S, x)$ is acting as a deck transformation pointwise on the values of the path $\widetilde{\eta}$ in $\widetilde{S}$. Note that $\widetilde{\eta}(0)$ and $\widetilde{\eta}(1)$ both descend to $x \in S$ since they both descend to $x_\circ \in \widetilde{S}/ N \subseteq \Sigma$ and $x_\circ$ descends to $x$. For $j \in \left\{0, 1 \right\}$, let $g_j: [0, 1] \rightarrow \widetilde{S}$ be the lift of $\gamma_1^{k_1}$ that starts at $\widetilde{\eta}(j)$. Then the path $\gamma_1^{k_1} \cdot \widetilde{\eta}$ is homotopic, relative to its endpoints, to the map $g_0^{-1}  \widetilde{\eta}  g_1$ where the concatenation is that of paths. This descends, as expected, to a loop in $\Sigma$ based at $\phi(x_\circ)$. Apply this to $\eta = a_{j, \ell} M \in \pi_1(\Sigma, x_\circ) \cong \pi_1(S, x)/ M$. Then we can view $\widetilde{\eta}$ as the path in $\widetilde{S}$ that goes from the constant path (in $S$) at $x$ to the path $\gamma_1^{\ell} \gamma_j^{-k_j} \gamma_1^{k_1 - \ell}$ (in $S$). Likewise, $g_0^{-1} \widetilde{\eta} g_1$ is the path that goes from $\gamma_1^{k_1}$ (the terminal point of $g_0$) to the path $\gamma_1^{\ell } \gamma_j^{-k_j} \gamma_1^{k_1 - (\ell-k_1)}$ (the terminal point of $g_1$). Since this does not start at the constant path, it does not descend to $\pi_1(\Sigma, x_\circ)$; instead, it descends to $\pi_1(\Sigma, \phi(x_\circ))$, as expected from the abstract theory telling us to expect a map of the form
$$\phi_*: \pi_1(\Sigma, x_\circ) \longrightarrow \pi_1(\Sigma, \phi(x_\circ)).$$
Use the change-of-basepoint map $\pi_1( \Sigma, \phi(x_\circ)) \rightarrow \pi_1(\Sigma, x_\circ)$ coming from the lift to $\Sigma$ of $\gamma_1^{-k_1}$ (this lift goes from $\phi(x_\circ)$ to $x_\circ$). At the universal cover level, this corresponds to left multiplying $g_0^{-1} \widetilde{\eta} g_1$ by $\gamma_1^{-k_1}$; this is a path from the constant path to $\gamma_1^{\ell - k_1} \gamma_j^{-k_j} \gamma_1^{k_1 - (\ell-k_1)}  = a_{j, \ell - k_1}$. From this, we see that the composition
$$\pi_1(\Sigma, x_\circ) \stackrel{\phi_*}{\longrightarrow} \pi_1(\Sigma, \phi(x_\circ)) \longrightarrow \pi_1(\Sigma, x_\circ)$$
acts by $a_{j, \ell} M \mapsto a_{j, \ell - k_1}M$. Taking $\alpha_{j, \ell} \in \pi_1(\Sigma, x_1)$ to be the image of $a_{j, \ell} M \in \pi_1(\Sigma, x_\circ)$ under the change-of-basepoint map coming from any path from $x_\circ$ to $x_1$, the theorem follows.

\bibliographystyle{alpha}

\end{document}